\documentclass[12pt]{report}
\usepackage{amssymb,amsmath,makeidx,verbatim}
\newcommand{\C}{\mathbb{C}}

\newcommand{\be}{\begin{enumerate}}
	\newcommand{\ee}{\end{enumerate}}
\newcommand{\bq}{\begin{eqnarray*}}
	\newcommand{\eq}{\end{eqnarray*}}

\begin{document}
	\newcommand{\disp}{\displaystyle}
	\thispagestyle{empty}
	\begin{center}
		\textsc{Fourier and Helgason Fourier transforms for Vector Bundle-valued
			Differential Forms on Homogeneous Spaces.\\}
\ \\
		\textsc{Olufemi O. Oyadare}\\
		Department of Mathematics,\\
		University of Ibadan,\\
		Ibadan, NIGERIA.\\
		\text{E-mail: \textit{femi$\_$oya@yahoo.com}}\\
	\end{center}
	\begin{quote}
		{\bf Abstract.} {\it The major aim of studying harmonic analysis on any group or symmetric space is to understand enough of its structure to be able to prove the fundamental theorem on it via the Fourier transform map. To this end we employ the perspective of the functional equation
			satisfied by the classical Fourier transform to derive the Helgason
			Fourier transform given as the map $\Omega^{l}(G/K,W)\longrightarrow\Omega^{k}(G/K\times G/P,V[\chi]):f\longmapsto \widehat{f}:G/K\times G/P\mapsto V[\chi]:(x,b)\longmapsto\widehat{f}(x,b)$ (for $W-$valued differential forms $f\in \Omega^{l}(G/K,W)$) as the $G-$ invariant vector bundle-valued differential form $\widehat{f}$ on the product
			space $G/K\times G/P$ whose image under the vector bundle-valued
			Poisson transform is the fibre convolution-integral $\varphi^{U^{\sigma,\nu}}_{\tau,l,k}* f$ on $G/K,$ where $\varphi^{U^{\sigma,\nu}}_{\tau,l,k}$ is the $W-$valued $\tau-$spherical $l-$form on $G/K.$ Explicitly, we prove that $$\widehat{f}_{l,k,\varepsilon(\lambda)}(x,b)=({\bf C_{o}(\lambda)}^{-1}\circ\beta^{V}(\lambda))\circ(\int_{G/K}\varphi^{U^{\sigma\nu},t}_{\lambda,l,k}\wedge\pi^{*}_{K}f)(x),$$ where $b\in G/P$ is a consequence of the boundary map $\beta^{V}(\lambda),$ ${\bf C_{o}(\lambda)}$ is the vector bundle-valued Harish-Chandra $c-$function and for some $\lambda-$linear relation, $\varepsilon(\lambda).$ The Fourier transform is found to be the map $\Omega^{l}(G/K,W)\longrightarrow\Omega^{k}(G/K\times G/P,W)$ $:f\mapsto f^{\triangle}:$ $G/P\times G/K\longrightarrow W$ $:(b,x)\longmapsto f^{\triangle}(b,x)$ and is then established to be explicitly given as $f^{\triangle}_{l,k,\upsilon(\lambda)}(b,x)=$ $$\int_{G/P}\phi_{k,l,\lambda}\wedge\pi^{*}_{P}(({\bf C_{o}(\lambda)}^{-1}\circ\beta^{V}(\lambda))\circ(\int_{G/K}\varphi^{U^{\sigma\nu},t}_{\lambda,l,k}\wedge\pi^{*}_{K}f)(x)),$$ where $\upsilon(\lambda)$ is some $\lambda-$linear relation. Our approach develops and distinguishes between the general theory of Fourier and Helgason Fourier transforms for the separate cases of invertible and non-invertible vector bundle-valued Poisson transforms on differential forms, in close relation with the calculus of differential forms of $G-$equivariant differential operators and the Bernstein-Gelfand-Gelfand (BGG-) sequences of strongly invariant differential operators on sections of the induced homology bundles, $\mathcal{H}_{k}(G/P,V).$}
	\end{quote}
	$\overline{2020\; \textmd{Mathematics}}$ Subject Classification: $42A38, \;\; 58A10, \;\; 53C30$
	\ \\
	\indent {\bf \S1. Introduction.}
	
	\indent This paper considers the general notion of {\it Fourier transform} as expounded in $[8.]$ and its direct relationship with the {\it Poisson transform} and the {\it Helgason Fourier transform} on non-compact symmetric spaces, and on the corresponding homogeneous vector bundles over such spaces and its compact dual. This approach opens up the contribution of homological algebra to the theory of these transforms and creates a general method for their construction in any function space of vector bundle-valued differential forms, as being sought in $[6.].$
	
	\indent We shall refer to a topological space $S$ as being {\it sufficiently structured} if the notions of {\it differentiability, integrability, spherical functions} and {\it convolution} are all well-defined on $S.$ By an arbitrary {\it Fourier transform} on $S$ we shall mean the map $$f\longmapsto f^{\triangle}(\lambda,s),$$ $s\in S$ and some $\lambda$ indexing the spherical functions $\varphi_{\lambda}$ on $S,$ in which $$f\longmapsto f^{\triangle}(\lambda,s):=(f*_{_{_{S}}}\varphi_{\lambda})(s),$$ $f\in C^{\infty}_{c}(S)$ and $*_{_{_{S}}}$ is the convolution on $S.$ The Euclidean space $\mathbb{R}^{n},$ every non-compact Riemannian symmetric space $S=G/K,$ every homogeneous vector bundle over $G/K$ and the homogeneous vector bundle over the compact Riemannian symmetric space $U/K$ dual to $G/K$ are all known examples of sufficiently structured topological spaces with corresponding Fourier transforms expressible as double integrals.
	
	\indent Indeed, for $S=\mathbb{R}^{n},$ we have $$f^{\triangle}(\lambda,\omega)=f^{\triangle}(\lambda\omega)=\int_{\mathbb{R}^{n}}f(x)e^{-i\langle x,u\rangle}dx=\int_{\mathbb{R}^{n}}\int_{(x,\omega)=p}e^{-ip\lambda}f(x)dm(x)dp$$ (where $dm$ is the Euclidean measure on the hyperplane $\xi=(x,\omega)=p.$ See $[5.],$ p. $10-11$); for a non-compact Riemannian symmetric space $S=G/K$ we have that $$f^{\triangle}(\lambda,x)=\int_{K/M}\int_{G/K}e^{(i\lambda+\rho)(A(x,b))}e^{(-i\lambda+\rho)(A(y,b))}f(y)dydb\;\;[5.];$$ for homogeneous vector bundles over a Riemannian symmetric space $G/K$ we also have that $$f^{\triangle}(\lambda,sx)=\frac{d_{\tau}}{d_{\bar{\sigma}'}}\int_{K}\int_{G}F^{i\lambda-\rho}(x^{-1}k)T_{\bar{\sigma}'}F^{i\bar{\lambda}-\rho}(y^{-1}k)^{*}f(y)dydk\;\;[1.]$$ and for the homogeneous vector bundle over the compact Riemannian symmetric space $U/K$ dual to $G/K,$ we have that $$f^{\triangle}(\mu,u)=\frac{d_{\tau}}{d_{\sigma}}\int_{K}\int_{U_{o}}e^{-\mu(A(k^{-1}u))}\tau(\mathbf{k}(k^{-1}u)^{-1})T_{\sigma,\mu}e^{(\mu+2\rho)(A(k^{-1}u_{o}))}f(u_{o})du_{o}dk$$ (with $supp\; f\subseteq U_{o};\; [2.]$).
	
	\indent Each of the above expressions for $f^{\triangle}(\lambda,s)$ is a double integral, which suggests writing it as a composition of two integral transforms. Indeed, a(n {\it incomplete}) transform in the first two cases of $S=\mathbb{R}^{n}$ (written as $f\mapsto \widehat{f}(\omega,p):=\int_{(x,\omega)=p}f(x)dm(x)$) and $S=G/K$ (written as $f\mapsto \widehat{f}(\lambda,b):=\int_{G/K}f(x)e^{(-i\lambda+\rho)(A(x,b))}dx$) leads to writing their earlier Fourier transforms (both $f^{\triangle}(\lambda\omega)$ and $f^{\triangle}(\lambda,x)$) as $$f^{\triangle}(\lambda\omega)=\int_{\mathbb{R}}\widehat{f}(\omega,p)e^{-ip\lambda}dp$$ and $$f^{\triangle}(\lambda,x)=\int_{K/M}e^{(i\lambda+\rho)(A(s,b))}\widehat{f}(\lambda,b)db,$$ respectively, from which the Poisson transforms $P_{\lambda}:F\longmapsto(P_{\lambda}F)(s),$ $s\in S,$ are seen as $$(P_{\lambda}F)(s)=\int_{\mathbb{R}}F(s,p)e^{-ip\lambda}dp$$ and $$(P_{\lambda}F)(s)=\int_{K/M}e^{(i\lambda+\rho)(A(s,b))}F(b)db,$$ respectively.
	
	\indent The present paper exploits the above intimate relationship, here seen as a {\it functional equation} of the type $$(f*_{_{_{S}}}\varphi_{\upsilon(\lambda)})(s)=P_{\lambda}(\widehat{f}(\varepsilon(\lambda),b)),$$ for some $\lambda-$linear relations $\varepsilon(\lambda)$ and $\upsilon(\lambda),$ relating the Fourier transform ($(f*_{_{_{S}}}\varphi_{\upsilon(\lambda)})(s)$), the Helgason Fourier transform ($\widehat{f}(\varepsilon(\lambda),b)$) and the Poisson transform $(P_{\lambda}),$ to start a discussion of their general theory, derive their expressions and deduce their properties in the context of vector bundle-valued differential forms. Here $\varepsilon(\lambda)=-\lambda$ for $G/K$ in the line bundle case ($[1.],$ p. $264$), $\varepsilon(\lambda)=-\lambda-i\rho$ for $G/K$ in the vector bundle case ($[1.],$ p. $275$) and $\varepsilon(\lambda)\in\Lambda(\tau):=\cup_{\sigma\in\widehat{M}(\tau)}\Lambda_{\sigma}(\tau)$ which is the set of highest restricted weights of all $U-$types in $\widehat{U}(\tau)$ for $U/K$ in the vector bundle case ($[2.],$ p. $102$).
	
	\indent The paper is organised as follows. The next section discusses the general theory of the Poisson transform for vector bundle-valued differential forms, while \S3. and \S4. give our main contribution. This includes a general theory for the Fourier and Helgason Fourier transforms and (for the special case of an invertible Poisson transform) the properties and explicit expressions for both of the Fourier and Helgason Fourier transforms.
	
	\ \\
	\indent {\bf \S2. The Poisson transform.}
	
	\indent Let $G$ be a connected semisimple Lie group with finite center, we denote its Lie algebra by $\mathfrak{g}$
	whose \textit{Cartan decomposition} is given as $\mathfrak{g} = \mathfrak{t}\oplus\mathfrak{p}.$ Denote by $\theta$ the \textit{Cartan involution} on $\mathfrak{g}$ whose collection of fixed points is $\mathfrak{t}.$
	We also denote by $K$ the analytic subgroup of $G$ with Lie
	algebra $\mathfrak{t}.$  $K$ is then a maximal compact subgroup of $G.$
	Choose a maximal abelian subspace  $\mathfrak{a}$ of $\mathfrak{p}$ with algebraic
	dual $\mathfrak{a}^*$ and set $A =\exp \mathfrak{a}.$  For every $\lambda \in \mathfrak{a}^*$ put
	$$\mathfrak{g}_{\lambda} = \{X \in \mathfrak{g}: [H, X] =
	\lambda(H)X, \forall  H \in \mathfrak{a}\},$$ and call $\lambda$ a restricted
	root of $(\mathfrak{g},\mathfrak{a})$ whenever $\mathfrak{g}_{\lambda}\neq\{0\}.$
	Denote by $\mathfrak{a}'$ the open subset of $\mathfrak{a}$
	where all restricted roots are $\neq 0,$ and call its connected
	components the \textit{Weyl chambers.}  Let $\mathfrak{a}^+$ be one of the Weyl
	chambers, define the restricted root $\lambda$ positive whenever it
	is positive on $\mathfrak{a}^+$ and denote by $\triangle^+$ the set of all
	restricted positive roots. Members of $\triangle^+$ which form a basis for $\triangle$ and can not be written as a linear combination of other members of $\triangle^+$ are called \textit{simple.} We then have the \textit{Iwasawa
		decomposition} $G = KAN$, where $N$ is the analytic subgroup of $G$
	corresponding to $$\mathfrak{n} = \sum_{\lambda \in \triangle^+} \mathfrak{g}_{\lambda},$$
	and the \textit{polar decomposition} $G = K\cdot
	cl(A^+)\cdot K,$ with $A^+ = \exp \mathfrak{a}^+,$ and $cl(A^{+})$ denoting the closure of $A^{+}.$
	
	If we set $$M = \{k \in K: Ad(k)H = H,\;H\in \mathfrak{a}\}$$ and $$M' = \{k
	\in K : Ad(k)\mathfrak{a} \subset \mathfrak{a}\}$$ and call them the
	\textit{centralizer} and \textit{normalizer} of $\mathfrak{a}$ in $K,$ respectively, then
	(i) $M$ and $M'$ are compact and have the same Lie algebra and
	(ii) the factor  $\mathfrak{w} = M'/M$ is a finite group called the \textit{Weyl
		group}. $\mathfrak{w}$ acts on $\mathfrak{a}^*_{\C}$ as a group of linear
	transformations by the requirement $$(s\lambda)(H) =
	\lambda(s^{-1}H),$$ $H \in \mathfrak{a}$, $s \in \mathfrak{w}$, $\lambda \in
	\mathfrak{a}^*_\mathbb{\C}$, the complexification of $\mathfrak{a}^*$.  We then have the
	\textit{Bruhat decomposition} $$G = \bigsqcup_{s\in \mathfrak{w}} B m_sB$$ where
	$B = MAN$ is a closed subgroup of $G$ (called the {\it minimal parabolic subgroup} of $G$) and $m_s \in M'$ is the
	representative of $s$ (i.e., $s = m_sM$). The Weyl group invariant members of a space shall be denoted by the superscript $^{\mathfrak{w}}.$ We fix an arbitrary {\it parabolic subgroup} $P$ of $G.$
	
	\indent Some of the most important functions on $G$ are the \textit{spherical
		functions} which we now discuss as follows.  A non-zero continuous
	function $\varphi$ on $G$ shall be called a \textit{(zonal) spherical
		function} whenever $\varphi(e)=1,$ $\varphi \in C(G//K):=\{g\in
	C(G)$: $g(k_1 x k_2) = g(x)$, $k_1,k_2 \in K$, $x \in G\}$ and $f*\varphi
	= (f*\varphi)(e)\cdot \varphi$ for every $f \in C_c(G//K),$ where $$(f \ast g)(x):=\int_{G}f(y)g(y^{-1}x)dy.$$ This
	leads to the existence of a homomorphism $\lambda :
	C_c(G//K)\rightarrow \C$ given as $\lambda(f) = (f*\varphi)(e)$.
	This definition is equivalent to the satisfaction of the functional relation $$\int_K\varphi(xky)dk = \varphi(x)\varphi(y),\;\;\;x,y\in G.$$
	
	\indent It has been shown by Harish-Chandra that spherical functions on $G$
	can be parametrized by members of $\mathfrak{a}^*_{\C}.$  Indeed every
	spherical function on $G$ is of the form $$\varphi_{\lambda}(x) = \int_Ke^{(i\lambda-p)H(xk)}dk,\; \lambda
	\in \mathfrak{a}^*_{\C},$$ with $\rho =
	\frac{1}{2}\sum_{\lambda\in\triangle^+} m_{\lambda}\cdot\lambda,$ where
	$m_{\lambda}=dim (\mathfrak{g}_\lambda),$ and that $\varphi_{\lambda} =
	\varphi_{\mu}$ iff $\lambda = s\mu$ for some $s \in \mathfrak{w}.$ Some of
	the well-known properties of spherical functions are $\varphi_{-\lambda}(x^{-1}) =
	\varphi_{\lambda}(x),$ $\varphi_{-\lambda}(x) =
	\bar{\varphi}_{\bar{\lambda}}(x),$ $\mid \varphi_{\lambda}(x) \mid\leq \varphi_{\Re\lambda}(x),$ $\mid \varphi_{\lambda}(x)\mid\leq \varphi_{i\Im\lambda}(x),$ $\varphi_{-i\rho}(x)=1,$ $\lambda \in \mathfrak{a}^*_{\C},$ while $\mid \varphi_{\lambda}(x) \mid\leq \varphi_{0}(x),\;\lambda \in i\mathfrak{a}^{*},\;x \in G.$ Also if $\Omega$ is the \textit{Casimir operator} on $G$ then
	$$\Omega\varphi_{\lambda} = -(\langle\lambda,\lambda\rangle +
	\langle \rho, \rho\rangle)\varphi_{\lambda},$$ where $\lambda \in
	\mathfrak{a}^*_{\C}$ and $\langle\lambda,\mu\rangle
	:=tr(adH_{\lambda} \ adH_{\mu})$ for elements $H_{\lambda}$, $H_{\mu}
	\in {\mathfrak{a}}.$ The elements $H_{\lambda}$, $H_{\mu}
	\in {\mathfrak{a}}$  are uniquely defined by the requirement that $\lambda
	(H)=tr(adH \ adH_{\lambda})$ and $\mu
	(H)=tr(adH \ adH_{\mu})$ for every $H \in {\mathfrak{a}}.$ Clearly $\Omega\varphi_0 = 0.$
	
	\indent A closer consideration of the properties of the elementary spherical functions (used in the definition of the scalar-valued Harish-Chandra Fourier transform) shows that $$\varphi_{-\lambda}(g)=\varphi_{\lambda}(g^{-1}),$$ for all $g\in G,$ $\lambda\in \mathfrak{a}^*_{\C}.$ This equality allows us to re-write the scalar-valued Fourier transform formula as $$\tilde{f}(\lambda)=\int_{G}f(g)\varphi_{-\lambda}(g)dg$$ 
	$$\;\;\;\;\;\;\;\;\;\;=\int_{G}f(g)\varphi_{\lambda}(g^{-1})dg$$ $$\;\;\;\;\;\;\;\;\;\;\;\;=\int_{G}f(g)\varphi_{\lambda}(eg^{-1})dg$$ $$\;\;=(f*_{_{_{G}}}\varphi_{\lambda})(e)\;\;$$
	$$=f^{\triangle}(\lambda,e),\;\;\;$$ where we have written $e$ for the identity element of $G$ and $*_{_{_{G}}}$ for the convolution of functions on $G.$ This computation shows that the scalar-valued Harish-Chandra Fourier transform is actually a convolution of every member of of the space $\mathcal{C}^{p}(G)$ with each of the elementary spherical functions on $G$ and that this convolution is then evaluated at $e.$ It would therefore be of general importance if we could consider a convolution of members of the space $\mathcal{C}^{p}(G)$ with each of the elementary spherical functions on $G$ (or with each of the {\it Eisenstein integral} on $G$) evaluated at any point of $G,$ which may not necessarily be at the identity element $e.$
	
	\indent Now let $\mathbb{W}_{\tau}$ denote a finite-dimensional $K-$representation with induced homogeneous vector bundle $W_{\tau}\longrightarrow G/K.$ We identify the space $\Gamma(W_{\tau})$ consisting of sections of the bundle with the space $C^{\infty}(G,\mathbb{W}_{\tau})^{K}$ consisting of smooth $K-$equivariant $\mathbb{W}_{\tau}-$valued functions on $G.$ Let $\mathbb{V}_{\tau}$ equally denote a finite-dimensional $P-$reprsentation with induced homogeneous vector bundle $V_{\tau}\longrightarrow G/P.$ We shall also identify the space $\Gamma(V_{\sigma})$ consisting of sections of the bundle $V_{\tau}\longrightarrow G/P$ with the space $C^{\infty}(G,\mathbb{V}_{\sigma})^{P}$ consisting of smooth $P-$equivariant $\mathbb{V}_{\sigma}-$valued functions on $G.$ If we include a $\rho-$shift into members of $C^{\infty}(G,\mathbb{V}_{\sigma})^{P}$ and define the representation space as $$C^{\infty}(G,\mathbb{V}_{\sigma})^{P,\rho}:=\{f\in C^{\infty}(G,\mathbb{V}_{\sigma})^{P}: f(gp)=e^{-\rho(H(p))}\sigma(p)^{-1}f(g)\},$$ then the space $C^{\infty}(G,\mathbb{V}_{\sigma})^{P,\rho}$ is isomorphic to the space of sections of the bundle $V_{\sigma}\otimes\mathcal{E}[\rho]\rightarrow G/P,$ in which $\mathcal{E}[\rho]$ is the density bundle over $G/P$ determined by $\rho\in\mathfrak{a}^{*}_{o}.$
	
	\indent The vector-valued Poisson transform associated to an $M-$equivariant linera map $\phi\in Hom_{M}(\mathbb{V}_{\tau},\mathbb{W}_{\tau})$ (called its {\it kernel}) shall be denoted by $$\Phi_{\phi}:C^{\infty}(G,\mathbb{V}_{\sigma})^{P,\rho}\longrightarrow C^{\infty}(G,\mathbb{W}_{\tau})^{K}$$ and is given as the integral $$\Phi_{\phi}(f)(gK)=\int_{K}(\tau(k)\circ\phi)(f(gk))dk,$$ for the normalized Haar measure $dk$ on $K.$ The $M-$invariance of the last integral makes it possible to change the last integration to integration over the Riemannian symmetric space $K/M.$ We now adapt the above definition of the Poisson transform to differential forms on $G/P.$
	
	\indent For the Lie algebra $\mathfrak{p}$ of the parabolic subgroup $P,$ let the corresponding $|k|-$grading of $\mathfrak{g}$ be given as $$\mathfrak{g}=\mathfrak{g}_{-k}\oplus\cdots\oplus\mathfrak{g}_{k}.$$ With $M=K\cap P,$ let the bundles associated to $\mathbb{W}$ (respectively, $\mathbb{V}$) over $G/M,$ be denoted as $W_{M}$ (respectively, $V_{M}$). For a $V-$valued $k-$form $\alpha$ on $G/P$ we can form the pullback of $\alpha$ along the projection $\pi_{P}:G/M\rightarrow G/P$ to arrive at a $\pi^{*}_{P}V-$valued $(0,k)-$form on $G/M.$ In particular we shall consider a $V[\chi]:=V\otimes\mathcal{E}[\chi]-$valued $k-$form $\alpha$ on $G/P,$ where $\mathcal{E}[\chi]$ denotes the density bundle corresponding to the character $\chi:P\rightarrow\mathbb{R}_{+}.$

\indent The space of all such $k-$forms would be denoted as $\Omega^{k}(G/P,V[\chi]).$ Now for any $Hom(V_{M},W_{M})-$valued $(l,n-k)-$form $\phi_{k,l}$ on $G/K\times G/P$ (in which $\dim(G/P)=n$) the wedge product $$\phi_{k,l}\wedge\pi^{*}_{P}\alpha$$ could be interpreted as a $W_{M}-$valued $(l,n)-$form on $G/K\times G/P$ by applying $\phi_{k,l}$ to $\pi^{*}_{P}\alpha$ after the insertion of vector fields on $G/K\times G/P.$ Hence, a fibre integral of $\phi_{k,l}\wedge\pi^{*}_{P}\alpha$ along the compact manifold $G/P$ would produce a $W-$valued $l-$form on $G/K,$ whose space is denoted as $\Omega^{l}(G/K,W).$ The foregoing argument produces a fibre-integral map from $\Omega^{k}(G/P,V[\chi])$ into $\Omega^{l}(G/K,W)$ here defined.
	
	\indent {\bf 2.1 Definition.} ($[4.]$) Let $\phi_{k,l}\in\Omega^{l,n-k}(G/K\times G/P,Hom(V_{M},W_{M})).$ The {\it Poisson transform} corresponding to $\phi_{k,l}$ written as $$\Phi_{\phi_{k,l}}:\Omega^{k}(G/P,V[\chi])\rightarrow\Omega^{l}(G/K,W)$$ is given by the fibre integral $$\Phi_{\phi_{k,l}}(\alpha)=\int_{G/P}\phi_{k,l}\wedge\pi^{*}_{P}\alpha.\;\;\;\Box$$
	
	\indent The notation for the Poisson transform corresponding to the $(l,n-k)-$form $\phi_{k,l}\in\Omega^{l,n-k}(G/K\times G/P,Hom(V_{M},W_{M}))$  given above shall be henceforth written simply as $\Phi_{k,l}.$ This construction answers the abstract homological approach to Poisson transform being sought in $[6].$ We shall refer to $\phi_{k,l}$ as the {\it kernel} of $\Phi_{k,l}$ and it can be shown that $\Phi_{k,l}$ is $G-$equivariant (that is, $g\cdot\Phi_{k,l}(\alpha)=\Phi_{k,l}(g\cdot\alpha),\;g\in G$) iff $\phi_{k,l}$ is $G-$invariant (that is, $g\cdot\phi_{k,l}=\phi_{k,l},\;g\in G$).
	
	\indent There is a bijective correspondence between the vector space of Poisson transforms $$\Phi_{k,l}:\Omega^{k}(G/P,V[\chi])\rightarrow\Omega^{l}(G/K,W)$$ and the vector space of all the $M-$invariant elements in the finite-dimensional representation $\Lambda^{l,n-k}(\mathfrak{g}/\mathfrak{m})^{*}\otimes L(\mathbb{V},\mathbb{W}).$ It is possible to also have a geometric description for the fibre-integral for the Poisson transform as contained in $[4],$ p. $48,$ which gives a more general outlook to the construction in $[3.].$
	
	\indent The above construction of the Poisson transform lends itself to the calculus of differential forms. Indeed, let $d^{W_{M}}$ and $d^{V_{M}}$ denote {\it covariant exterior derivatives} on $W_{M}$ and $V_{M},$ respectively, induced by the respective tractor connections which splits into $d^{W_{M}}=d^{W_{M}}_{K}+d^{W_{M}}_{P}$ and $d^{V_{M}}=d^{V_{M}}_{K}+d^{V_{M}}_{P}.$ Let $$*:\Omega^{k}(G/K,V)\rightarrow\Omega^{\dim(G/K)-k}(G/K,V)$$ denote the {\it Hodge star operator} given by the formula $$\alpha\wedge_{h}*\beta=\langle\alpha,\beta\rangle vol,$$ for $\alpha,\beta\in\Omega^{k}(G/K,V)$ and obtain the {\it $P-$codifferential tensorial operator} $$*_{K}:\Lambda^{p,q}T^{*}(G/K\times G/P)\otimes Z\rightarrow\Omega^{\dim(G/K)-p,q}T^{*}(G/K\times G/P)\otimes Z$$ defined on $\alpha\wedge\beta$ (where we have chosen $\alpha\in\Lambda^{p,0}T^{*}_{x}(G/M)\otimes(V_{M})_{x}$ and $\beta\in\Lambda^{0,q}T^{*}_{x}(G/K/\times G/P)\otimes(W_{M})^{*}_{x},$ $x\in G/M$) as $$*_{K}(\alpha\wedge\beta)=(*\alpha)\wedge\beta.$$ We shall also consider the $P-$codifferential of the {\it Kostant codifferential} $\partial^{*}$ given as $$\partial_{P}:\Lambda^{p,q}T^{*}(G/M)\otimes Z\rightarrow\Omega^{p,q-1}T^{*}(G/M)\otimes Z$$ and defined as $$\partial_{P}(\alpha\wedge\beta)=(-1)^{p}\alpha\wedge\partial^{*}\beta,$$ for $\alpha\in\Lambda^{p,0}T^{*}_{x}(G/M)\otimes(V_{M})_{x}$ and $\beta\in\Lambda^{0,q}T^{*}_{x}(G/M)\otimes(W_{M})^{*}_{x},$ $x\in G/M.$ The covariant codifferential $$\delta^{W}:\Omega^{l}(G/K,W)\rightarrow\Omega^{l-1}(G/K,W)$$ is defined by $$\delta^{W}\alpha:=(-1)^{k}*^{-1}d^{\theta}*\alpha,$$ where $*$ is the Hodge star operator and $d^{\theta}$ is the covariant exterior derivative on $W_{M}-$valued differential forms on $G/K$ induced by the {\it $\theta-$twisted tractor connection.}
	
	\indent It can be proved (see $[4.],$ p. $52$) that for any Poisson transform $\Phi_{k,l},$ the compositions $d^{W_{M}}\circ\Phi_{k,l},$ $*\circ\Phi_{k,l}$ and $\delta^{W}\circ\Phi_{k,l}$ are all Poisson transforms with corresponding kernels $d^{W_{M}}_{K}\phi_{k,l},$ $*_{K}\phi_{k,l}$ and $\delta^{W}_{K}\phi_{k,l},$ respectively, while the compositions $\Phi_{k,l}\circ d^{V}$ and $\Phi_{k,l}\circ\partial^{*}$ are equally Poisson transforms with corresponding kernels given as $(-1)^{n-k+l+1}d^{W_{M}}_{P}\phi_{k,l}$ and $(-1)^{n-k+l}\partial_{P}\phi_{k,l},$ respectively.
	
	\ \\
	\indent {\bf \S3. The Fourier and Helgason Fourier transforms.}
	
	\indent Consider now a $W-$valued $l-$form $f$ on the Riemannian symmetric space $G/K,$ from which we form the pullback $\pi^{}_{K}f$ of $f$ along the projection $\pi_{K}:G/M\rightarrow G/K$ to arrive at a $\pi^{*}_{K}-$valued $(l,0)-$form on $G/M.$ If $\sigma\in\widehat{M}$ so that $\sigma$ occurs with multiplicity $m_{\sigma}>0$ in $\tau_{|_{_{M}}}$ and $\nu:\mathfrak{a}\rightarrow\mathbb{C}$ is a linear function on $\mathfrak{a},$ we set $$U^{\sigma\nu}=ind^{G}_{P}(\sigma\otimes\exp(\nu+\rho)\otimes\bf{1})$$ as the nonunitary (minimal) principal series of $G.$ We fix $\sigma$ and $\nu$ such that $\tau\subseteq U^{\sigma\nu}_{|_{K}}$ and denote the corresponding $Hom(V_{M},W_{M})-$valued $(m-l,k)-$ $\tau-$spherical form on $G/K\times G/P$ (in which $m=\dim(G/K)$) by $\varphi^{U^{\sigma\nu}}_{\tau,l,k}$ whose translation is also written as $\varphi^{U^{\sigma\nu},t}_{\tau,l,k}.$ See $[1.]$ for the classical motivation from harmonic analysis.
	
	\indent The wedge product $$\varphi^{U^{\sigma\nu},t}_{\tau,l,k}\wedge\pi^{*}_{K}f$$ is then interpreted as a $V_{M}-$valued $(m,k)-$form on $G/K\times G/P,$ by applying $\varphi^{U^{\sigma\nu},t}_{\tau,l,k}$ to $\pi^{*}_{K}f$ after the insertion of vector fields on $G/K\times G/P.$ Hence, a fibre-integral of $$\varphi^{U^{\sigma\nu},t}_{\tau,l,k}\wedge\pi^{*}_{K}f$$ along the oriented manifold $G/K$ would produce a $V-$valued $k-$form on $G/P.$ The above creates a fibre-integral map from $\Omega^{l}(G/K,W)$ into $\Omega^{k}(G/P,V)$ which leads to the following definition.\\
	
	\indent {\bf 3.1 Definition.} Let $\varphi^{U^{\sigma\nu}}_{\tau,l,k}$  be as above and let $f\in\Omega^{l}(G/K,W).$ The {\it convolution} $\varphi^{U^{\sigma\nu}}_{\tau,l,k}*f$ is defined as the $V-$valued $k-$form on $G/K$ given by the fibre-integral $$\varphi^{U^{\sigma\nu}}_{\tau,l,k}*f=\int_{G/K}\varphi^{U^{\sigma\nu},t}_{\tau,l,k}\wedge\pi^{*}_{K}f.\;\;\;\;\Box$$
	
	\indent The above convolution helps to define the {\it Helgason Fourier transform} as suggested by our observations in the Introduction.\\
	
	\indent {\bf 3.2 Definition.} Let $f\in\Omega^{l}(G/K,W).$ The {\it Helgason Fourier transform} of $f$ is the map $$\Omega^{l}(G/K,W)\longrightarrow\Omega^{k}(G/K\times G/P,V[\chi])$$ $$\;\;\;\;\;\;\;\;\;\;\;\;\;\;\;\;\;\;:f\longmapsto \widehat{f}:G/K\times G/P\mapsto V[\chi]$$ $$\;\;\;\;\;\;\;\;\;\;\;\;\;\;\;\;\;\;\;\;\;\;\;\;\;\;\;\;\;\;\;\;:(x,b)\longmapsto\widehat{f}_{l,k}(x,b)\;\;\;$$ for which $$(\varphi^{U^{\sigma\nu}}_{\tau,l,k}*f)(x)=\Phi_{k,l}(\widehat{f}_{l,k}(x,b)),$$ where $x\in G/K$, with the $V[\chi]-$valued form $\widehat{f}_{l,k}(x,b)$ considered as a function of $b\in G/P$  and $\Phi_{k,l}$ is the Poisson transform (corresponding to some kernel $\phi_{k,l}\in\Omega^{l,n-k}(G/K\times G/P,Hom(V_{M},W_{M}))$)$.\;\;\;\Box$
	
	\indent It should be noted here that the Poisson transform in the last definition is a fibre integral of $\widehat{f}_{l,k}(x,b)$ over $b\in G/P,$ thus making the convolution $\varphi^{U^{\sigma\nu}}_{\tau,l,k}*f$ a $W-$valued $l-$form on $x\in G/K$ and making the Helgason Fourier transform $f\mapsto\widehat{f}_{l,k}$ as being defined on $(x,b)\in G/K\times G/P.$
	
	\indent In order words, $\alpha\in\Omega^{k}(G/K\times G/P,V[\chi])$ is said to be the Helgason Fourier transform of $f\in\Omega^{l}(G/K,W)$ (written as $\alpha=\widehat{f}_{l,k}$) if the Poisson-transform image, $\Phi_{k,l}(\alpha),$ of $\alpha$ is exactly the convolution $\varphi^{U^{\sigma\nu}}_{\tau,l,k}*f,$ for some $Hom(V_{M},W_{M})-$valued $(m-l,k)-$ $\tau-$spherical form, $\varphi^{U^{\sigma\nu}}_{\tau,l,k},$ on $G/K\times G/P.$ Hence, the Helgason Fourier transform $\widehat{f}_{l,k}(x,b)$ of $f$ is the $V[\chi]-$valued map defined on all $(x,b)\in G/K\times G/P$ which, when considered as a function on only all $b\in G/P,$ has its Poisson transform to be exactly the fibre-integral convolution $(\varphi^{U^{\sigma\nu}}_{\tau,l,k}*f)(x).$ This is a far reaching characterization of the notion of Helgason Fourier transform.
	
	\indent Even though Definition $3.2$ implies that the Poisson transform and the Helgason Fourier transform are inter-dependent on each other (and hence, together dependent on the kernel $\phi_{k,l}$), the existence of the Helgason Fourier transform is not dependent on the invertibility of the Poisson transform. It appears however that an explicit fibre-integral expression of the Helgason Fourier transform would require the use of the inverse, $\Phi^{-1}_{k,l},$ so that we could have that $$\widehat{f}_{l,k}(x,b)=\Phi^{-1}_{k,l}((\varphi^{U^{\sigma\nu}}_{\tau,l,k}*f)(x))$$ if and when $\Phi^{-1}_{k,l}$ exists. Though we shall in the next section discuss the nature of $\Phi^{-1}_{k,l},$ Definition $3.2$ suggests that the Helgason Fourier transform may be completely defined and adequately studied with or without the invertibility of the Poisson transform.
	
	\indent Definition $3.1$ may also be considered as a genuine transform on members $f\in\Omega^{l}(G/K,W).$
	
	\indent {\bf 3.3 Definition.} Let $f\in\Omega^{l}(G/K,W).$ The {\it Fourier transform,} $f^{\triangle}_{l,k}(b,x),$ $b\in G/P,$ $x\in G/K,$ of $f$ is the map $$\Omega^{l}(G/K,W)\rightarrow\Omega^{k}(G/K\times G/P,W)\;\;\;\;\;\;\;\;\;\;\;\;\;$$  $$\;\;\;f\mapsto f^{\triangle}_{l,k}(\cdot,\cdot):G/P\rightarrow W$$ $$\;\;\;\;\;\;\;\;\;\;\;\;\;\;\;\;\;\;\;\;\;\;\;\;\;\;\;\;\;\;\;\;\;\;\;\;\;\;\;\;\;\;\;\;\;\;\;\;\;\;\;\;\;:b\mapsto f^{\triangle}_{l,k}(b,x),\;\;\mbox{with}\;\;x\in G/K,$$ for which $$f^{\triangle}_{l,k}(b,x):=(\varphi^{U^{\sigma\nu}}_{\tau,l,k}*f)(x).\;\;\;\Box$$
	
	\indent Clearly, $f^{\triangle}_{l,k}(b,x)=\Phi_{k,l}(\widehat{f}_{l,k}(x,b)).$ The Fourier transform is, in its full essence, the map $$\Omega^{l}(G/K,W)\longrightarrow\Omega^{k}(G/K\times G/P,W).$$ This is the wholistic generalization of the classical examples of the {\it Fourier transform} to vector bundle-valued differential form on $G/K$ as earlier posited in the Introduction, whose definition was there shown to be interwoven with those of the Poisson transforms. It subsumes the {\it Joint-Eigenspace Fourier transform} studied in $[8.],\;[9.]$ and $[10.].$ This definition presents the Fourier transform as the Poisson-transform-completion of the Helgason Fourier transform.
	
	\indent We may also refer to $f^{\triangle}_{l,k}(b,x)$ as the Fourier transform with respect to the Poisson kernel $\phi_{k,l}\in\Omega^{l,n-k}(G/K\times G/P,Hom(V_{M},W_{M}))$ in which we have written $f^{\triangle}_{l,k}(b,x)$ for $f^{\triangle}_{\phi_{k,l}}(b,x).$
	
	\indent It is noted here that, as against the Poisson transform $$\Phi_{k,l}:\Omega^{k}(G/P,V[\chi])\rightarrow\Omega^{l}(G/K,W)$$ (as a map from the {\it Furstenberg boundary} $G/P$ to all of $G/K$) and the Helgason Fourier transform $$\Omega^{l}(G/K,W)\longrightarrow\Omega^{k}(G/K\times G/P,V[\chi])$$ (being a map from all of $G/K$ to $G/K\times G/P$), the Fourier transform is a map $$\Omega^{l}(G/K,W)\rightarrow\Omega^{k}(G/K\times G/P,W)$$ among the $W-$valued differential forms of different orders. This is reminiscence of the self-duality of the classical Fourier transform on $\mathbb{R}^{n}.$ Furthermore, the bijective correspondence between the space of Poisson transforms and the vector space of $M-$invariant elements in the finite-dimensional representation $\Lambda^{l,n-k}(\mathfrak{g}/\mathfrak{m})^{*}\otimes L(\mathbb{V},\mathbb{W})$ transfers to the vector space of all Fourier transforms via Definition $3.2.$
	
	\indent It may sometimes be convenient to write $f^{\triangle}_{l,k}$ simply as $\Phi_{k,l}\circ\widehat{f}_{l,k},$ since the space $\Omega^{k}(G/K\times G/P,W)$ may be viewed as $\Omega^{k}(G/P\times G/K,W),$ so that whenever the Poisson transform is inverted, we would have that $$\widehat{f}_{l,k}=\Phi^{-1}_{k,l}\circ f^{\triangle}_{l,k}.$$
	
	\indent Our objective is to study the Fourier and Helgason Fourier transform in the general setting of a vector bundle-valued differential forms, in this section for any Poisson transform and later in the next section, for invertible Poisson transforms. We however start with the calculus of Fourier transform.
	
	\indent {\bf 3.4 Theorem.} Let $Z$ be the homogeneous vector bundle over $G/M$ corresponding to the $M-$representation $L(\mathbb{V},\mathbb{W})$ and let $\phi_{k,l}$ be a $Z-$valued Poisson transform kernel of degree $(l,n-k)$ with the corresponding Poisson transform $$\Phi_{k,l}:\Omega^{k}(G/P,V[\chi])\rightarrow\Omega^{l}(G/K,W).$$
	
	\indent $(i)$ The compositions $d^{W}\circ f^{\triangle}_{l,k},$ $*\circ f^{\triangle}_{l,k}$ and $\delta^{W}\circ f^{\triangle}_{l,k}$ are all Fourier transforms of $f\in\Omega^{l}(G/K,W)$ with respect to the corresponding Poisson kernels $d_{K}\phi_{k,l},$ $*_{K}\phi_{k,l}$ and $\delta_{K}\phi_{k,l},$ respectively.
	
	\indent $(ii)$ The compositions $f^{\triangle}_{l,k}\circ d^{V}$ and $f^{\triangle}_{l,k}\circ\partial^{*}$ are Fourier transforms of $f\in\Omega^{l}(G/K,W)$ with respect to the corresponding Poisson kernels given as $(-1)^{n-k+l+1}d_{P}\phi_{k,l}$ and $(-1)^{n-k+l}\partial_{P}\phi_{k,l},$ respectively.
	
	\indent {\bf Proof.} We give the proof for the case of the composition $d^{W}\circ f^{\triangle}_{l,k},$ since the remaining ones are established in much the same manner. This may be compared with Proposition $2.2.3$ of $[4].$ Indeed, $d^{W}\circ f^{\triangle}_{(d_{K}\phi_{k,l})}(b,x)=d^{W}\circ f^{\triangle}_{l,k}(b,x)=d^{W}(\Phi_{k,l}(\widehat{f}_{l,k}(x,b)))$ $$=d^{W}(\int_{G/P}\phi_{k,l}\wedge\pi^{*}_{P}(\widehat{f}_{l,k}(x,b)))=\int_{G/P}(d_{K}\phi_{k,l})\wedge\pi^{*}_{P}(\widehat{f}_{l,k}(x,b))$$ $$=\Phi_{(d_{K}\phi_{k,l})}(\widehat{f}_{l,k}(x,b))=f^{\triangle}_{(d_{K}\phi_{k,l})}(b,x).$$ The remaining compositions are established in the same manner$.\;\;\Box$
	
	\indent The Fourier transform $f^{\triangle}_{l,k}(b,x)$ is intimately related (by Definitions $3.2$ and $3.3$) to the Helgason Fourier transform $\widehat{f}_{l,k}(x,b)$ of every $f\in\Omega^{l}(G/K,W).$ This is the major reason some authors (of the Helgason school of thought) could not distinguish between the two transforms. One of such intimate relationships is the following results.
	
	\indent {\bf 3.5 Proposition.} The set of eigenvalues $\mathfrak{E}_{D_{C}}(f^{\triangle}_{l,k})$ of the Fourier transform $f\rightarrow f^{\triangle}_{l,k}$ is exactly the same as the set of eigenvalues $\mathfrak{E}_{D_{C}}(\widehat{f}_{l,k})$ of the Helgason Fourier transform $f\rightarrow \widehat{f}_{l,k},$ under the $G-$equivariant differential operator $D_{C}$ of the Casinir operator.
	
	\indent {\bf Proof.} Let $\lambda_{D_{C}}(f^{\triangle}_{l,k})\in\mathfrak{E}_{D_{C}}(f^{\triangle}_{l,k})$ and $\lambda_{D_{C}}(\widehat{f}_{l,k})\in\mathfrak{E}_{D_{C}}(\widehat{f}_{l,k}).$ We then have that $$\lambda_{D_{C}}(f^{\triangle}_{l,k})\cdot f^{\triangle}_{l,k}(b,x)=D_{C}f^{\triangle}_{l,k}(b,x)\;\;\;\;\;\;\;\;\;\;\;\;\;\;\;\;\;\;\;\;\;\;\;\;\;\;\;\;\;\;\;\;\;\;\;\;\;\;\;\;\;\;\;\;\;\;$$ $$\;\;\;\;\;\;\;\;\;\;\;\;\;\;\;\;\;\;\;\;\;\;\;\;\;=D_{C}(\Phi_{k,l}(\widehat{f}_{l,k}(x,b))),\mbox{(by Definition $3.3$)}$$ $$=\Phi_{k,l}(D_{C}(\widehat{f}_{l,k}(x,b))),\;\;\;\;$$ (due to the $G-$invariance of the Poisson transform) $$\;\;\;\;\;=\Phi_{k,l}(\lambda_{D_{C}}(\widehat{f}_{l,k})\cdot\widehat{f}_{l,k}(x,b))$$ $$\;=\lambda_{D_{C}}(\widehat{f}_{l,k})\cdot(\Phi_{k,l}(\widehat{f}_{l,k}))$$ $$\;\;\;\;\;\;\;\;\;\;\;\;\;\;\;\;\;\;\;\;\;\;\;\;\;\;\;\;\;=\lambda_{D_{C}}(\widehat{f}_{l,k})\cdot f^{\triangle}_{l,k}(b,x),\mbox{(by Definition $3.3$)}.$$ Hence, $\lambda_{D_{C}}(f^{\triangle}_{l,k})=\lambda_{D_{C}}(\widehat{f}_{l,k}).$ Thus, $$\mathfrak{E}_{D_{C}}(f^{\triangle}_{l,k})=\{\lambda_{D_{C}}(f^{\triangle}_{l,k})\in \mathbb{C}:D_{C}f^{\triangle}_{l,k}(b,x)=\lambda_{D_{C}}(f^{\triangle}_{l,k})\cdot f^{\triangle}_{l,k}(b,x)\}$$ $$\;\;\;\;\;\;\;\;\;\;\;\;\;\;\;=\{\lambda_{D_{C}}(\widehat{f}_{l,k})\in \mathbb{C}:D_{C}\widehat{f}_{l,k}(b,x)=\lambda_{D_{C}}(\widehat{f}_{l,k})\cdot \widehat{f}_{l,k}(b,x)\}$$ $$=\mathfrak{E}_{D_{C}}(\widehat{f}_{l,k}).\;\;\Box\;\;\;\;\;\;\;\;\;\;\;\;\;\;\;\;\;\;\;\;\;\;\;\;\;\;\;\;\;\;\;\;\;\;\;\;\;\;\;\;\;\;\;\;$$
	
	\indent What we have actually shown in the last result is the stronger assertion that a scalar is an eigenvalue of the Fourier transform under $D_{C}$ if and only if the {\it same scalar} is an eigenvalue of the Hegason Fourier transform under $D_{C}.$
	
	\indent It is clear that more structural properties of the Poisson transform would be required in order to understand more of the nature of the Fourier and Helgason Fourier transforms. A correspondence calculus between the {\it covariant Laplacian} and {\it the curved Laplacian} sets the stage for an application of the {\it Bernstein-Gelfand-Gelfand (BGG-)compatibility} to the theory.
	
	\indent {\bf 3.6 Lemma.} Let $\triangle^{W}$ denote the {\it covariant Laplace operator} (also called the {\it twisted Laplace-Beltrami} operator) on the space $\Omega^{l}(G/K,W)$ given by $\triangle^{W}=d^{W}\delta^{W}+\delta^{W}d^{W}$ and let $\Box^{R}$ denote the {\it curved Laplace operator} on $\Omega^{k}(G/P,V)$ given by $\Box^{R}=d^{V}\partial^{*}+\partial^{*}d^{V}.$ If $c_{\mathbb{W}}$ and $c_{\mathbb{V}}$ represent the Casimir eigenvalues of the Casimir induced differential operator $D_{C}$ of Proposition $3.5$ then, for every Poisson transform $\Phi_{k,l}:\Omega^{k}(G/P,V)\rightarrow\Omega^{l}(G/K,W),$ $\triangle^{W}$ acts on $\Phi_{k,l}(\alpha)$ as $c_{\mathbb{W}}-c_{\mathbb{V}}-2c_{\Box^{R}}$ if and only if $\Box^{R}$ acts on $\alpha$ as $c_{\Box^{R}}.$
	
	\indent {\bf Proof.} It is known ($[4.],\;p.\;54$) that the differential operator $D_{C}$ acts on $V-$valued differential forms on $G/P$ as $$D_{C}=2\Box^{R}+c_{\mathbb{V}},$$ while it acts on $W-$valued differential forms on $G/K$ as $$D_{C}=-\triangle^{W}+c_{\mathbb{W}}.$$ Since $D_{C}\Phi_{k,l}(\alpha)=\Phi_{k,l}(D_{C}\alpha),$ it follows that $$(-\triangle^{W}+c_{\mathbb{W}})\Phi_{k,l}(\alpha)=\Phi_{k,l}(2\Box^{R}+c_{\mathbb{V}}).$$ That is, $$-\triangle^{W}\Phi_{k,l}(\alpha)+c_{\mathbb{W}}\Phi_{k,l}(\alpha)=2\Phi_{k,l}(\Box^{R}\alpha)+c_{\mathbb{V}}\Phi_{k,l}(\alpha).$$ Therefore, $$\triangle^{W}\Phi_{k,l}(\alpha)=(c_{\mathbb{W}}-c_{\mathbb{V}})\Phi_{k,l}(\alpha)-2\Phi_{k,l}(\Box^{R}\alpha).$$ Now, $\triangle^{W}$ acts on $\Phi_{k,l}(\alpha)$ as $c_{\mathbb{W}}-c_{\mathbb{V}}-2c_{\Box^{R}}$ if and only if $$(c_{\mathbb{W}}-c_{\mathbb{V}}-2c_{\Box^{R}})\Phi_{k,l}(\alpha)=(c_{\mathbb{W}}-c_{\mathbb{V}})\Phi_{k,l}(\alpha)-2\Phi_{k,l}(\Box^{R}\alpha),$$ if and only if $$\Phi_{k,l}(\Box^{R}\alpha)=c_{\Box^{R}}\Phi_{k,l}(\alpha)=\Phi_{k,l}(c_{\Box^{R}}\alpha),$$ if and only if $$\Phi_{k,l}(\Box^{R}\alpha-c_{\Box^{R}}\alpha)=0\;\;\mbox{(due to linearity of $\Phi_{k,l}$)},$$ if and only if $$\Box^{R}\alpha=c_{\Box^{R}}\alpha,$$ due both to the fact that $\Phi_{k,l}$ is well-defined (i.e., $\alpha_{1}=\alpha_{2}\implies\Phi_{k,l}(\alpha_{1})=\Phi_{k,l}(\alpha_{2})$) and the non-degeneracy of $\Phi_{k,l}$ (i.e., $\Phi_{k,l}(\alpha_{1}-\alpha_{2})=0\implies\alpha_{1}-\alpha_{2}=0)$ from the non-degenracy of the wedge product in Definition $2.1.\;\Box$

\indent The exact form of the eigenvalue $c_{\Box^{R}}$ for $f^{\triangle}_{l,k}(b,x)$ may be inferred from the contributions of $d^{V}$ and $\partial^{*}$ (given in Theorem $3.4$) into $\Box^{R}=d^{V}\partial^{*}+\partial^{*}d^{V}.$ 
	
	\indent We now come to the important application of the {\it Bernstein-Gelfand-Gelfand (BGG) sequences} of strongly invariant differential operators on sections of the {\it induced homology bundles,} $\mathcal{H}_{k}(G/P,V),\;[4.].$ The geometric properties of the Fourier transform are deeply influenced by its satisfaction of the Bernstein-Gelfand-Gelfand (BGG-)compatibility defined below.
	
	\indent {\bf 3.7 Definitions.}\\
\indent $1.$ A {\it $W-$valued harmonic differential} on $G/K\times G/P$ is a differential form $\alpha\in\Omega^{k}(G/K\times G/P,W)$ in which $\Box^{R}\cdot\alpha=0.$ In this case we say that $\alpha$ is contained in the $W-$valued harmonic differential forms on $G/K\times G/P.$

\indent $2.$ The Fourier transform $f\rightarrow f^{\triangle}_{l,k}(b,x)$ (as given in Definition $3.3$) is said to be {\it Bernstein-Gelfand-Gelfand (BGG-)compatible} whenever it is contained in the $W-$valued harmonic differential forms on $G/K\times G/P.\;\;\Box$
	
	\indent The following lemma prepares the ground for a major result on BGG-compatible Fourier transforms.
	
	\indent {\bf 3.8 Lemma.} The Fourier transform $f\rightarrow f^{\triangle}_{l,k}(b,x)$ is BGG-compatible if and only if its corresponding Poisson transform $\alpha\mapsto\Phi_{k,l}(\alpha)$ is BGG-compatible.
	
	\indent {\bf Proof.} We recall the functional equation $$f^{\triangle}_{l,k}(b,x)=\Phi_{k,l}(\widehat{f}_{l,k}(x,b))$$ and note (from Theorem $2.2.3$ of $[4.]$) that {\it 'the image of the transform $\Phi$ is contained in the harmonic differential forms'} is equivalent to {\it 'the image of the transform $f^{\triangle}$ is contained in the harmonic differential forms.'} Hence, the BGG-compatibility of $f^{\triangle}_{l,k}(b,x)$ is equivalent to the BGG-compatibility of $\Phi_{k,l}(\alpha).\;\;\Box$
	
	\indent The culmination of the last three results is the following characterization of the BGG-compatible Helgason Fourier transforms.
	
	\indent {\bf 3.9 Theorem.} Every BGG-compatible Helgason Fourier transform $$f\rightarrow \widehat{f}_{l,k}(b,x)$$ has that its corresponding Fourier transform satisfies the differential equation $$\triangle^{W}\cdot f^{\triangle}_{l,k}(b,x)=(c_{\mathbb{W}}-c_{\mathbb{V}})f^{\triangle}_{l,k}(b,x).$$
	
	\indent {\bf Proof.} 

\indent We know from Theorem $3.6$ that $\triangle^{W}$ acts on $\Phi_{k,l}(\alpha)$ as $c_{\mathbb{W}}-c_{\mathbb{V}}-2c_{\Box^{R}}$ if and only if $\Box^{R}$ acts on $\alpha$ as $c_{\Box^{R}}.$ That is, $$\triangle^{W}\cdot\Phi_{k,l}(\alpha)=(c_{\mathbb{W}}-c_{\mathbb{V}}-2c_{\Box^{R}})\Phi_{k,l}(\alpha)\;\;\mbox{if and only if}\;\;\Box^{R}\cdot\alpha=c_{\Box^{R}}\alpha.$$ With $\alpha$ set as $\widehat{f}_{l,k}(b,x),$ and using $$0=\Box^{R}\cdot\widehat{f}_{l,k}(b,x)=c_{\Box^{R}}\widehat{f}_{l,k}(b,x) \;\;\;\mbox{(from the hypothesis),}$$ we hve that $c_{\Box^{R}}=0.$ Hence, $$\triangle^{W}\cdot f^{\triangle}_{l,k}(b,x)=(c_{\mathbb{W}}-c_{\mathbb{V}})f^{\triangle}_{l,k}(b,x).\;\Box$$

	\indent {\bf 3.10 Corollary.} In the special case that $\mathbb{W}=\mathbb{V}$ (so that we have that $\triangle^{W}=\triangle^{V}=\triangle$), the Fourier transform of Theorem $3.9$ satisfies the differential equation $$\triangle\cdot f^{\triangle}_{l,k}(b,x)=0.\;\;\Box$$
	
	\indent The special scalar-valued case of Corollary $3.10$ leads to the consideration of the joint-eigenspace of harmonic analysis (for both the Fourier transform $f^{\triangle}_{l,k}$ and the Helgason Fourier transform $\widehat{f}_{l,k}$ of $f$) by setting the differential operator $\triangle$ as $$\triangle=D-\Gamma(D)(i\lambda)$$ (for every $D\in {\bf D}(G/K)$) on a BGG-compatible $f^{\triangle}_{l,k}$ and $\widehat{f}_{l,k},$ in order to have the differential equations $$D\cdot f^{\triangle}_{l,k}=\Gamma(D)(i\lambda)f^{\triangle}_{l,k}\;\;\mbox{and}\;\;D\cdot \widehat{f}_{l,k}=\Gamma(D)(i\lambda)\widehat{f}_{l,k},$$ respectively. Indeed, the inversion and Plancherel formula, and Paiey-Wiener theorem for the Fourier transform of Definition $3.3$ (named the {\it Joint-Eigenspace Fourier transform}) of the scalar-valued $0-$differential forms on $G/K$ were the main results of $[8.],$ which would be established for all $\Omega^{l}(G/K,W)$ in a forthcoming paper.
	
	\indent {\bf 3.11 Remarks.}
	
	\indent $1.$ The known theory of the BGG-sequences was crafted for the parabolic geometries $G/P.$ Hence its direct application to the Fourier and Poisson transforms, but not to the Helgason Fourier transform. In order to discuss more geometric properties of the Helgason Fourier transform, there is the need to develop a corresponding theory of the BGG-sequences for the homogeneous space $G/K \times G/P$ and to then seek the image of the map $f\mapsto\widehat{f}_{l,k}$ among the $Hom(V_{M},W_{M})-$valued $(l,n-k)-$forms on $G/K \times G/P\eqsim G/M.$ This quest reduces to the study of $M-$invariant elements in the finite-dimensional representation $\Lambda^{l,n-k}(\mathfrak{g}/\mathfrak{m})^{*}\otimes L(\mathbb{V},\mathbb{W}).$
	
	\indent $2.$ In the meantime, it is convenient to refer to the Helgason Fourier transform as being BGG-compatible whenever its Poisson transform is BGG-compatible.
	
	\indent $3.$ The theory of the Fourier and Helgason Fourier transforms of this section is independent of the invertibility of the Poisson transform. Indeed, we solely depended on the functional equation $$f^{\triangle}_{l,k}(b,x)=\Phi_{k,l}(\widehat{f}_{l,k}(x,b))=\int_{G/P}\phi_{k,l}\wedge\pi^{*}_{P}(\widehat{f}_{l,k}(x,b))$$ of Definition $3.3$ to proceed this far. In particular, we have not given any explicit expression for each of the Fourier and Helgason Fourier transforms, despite the fact that the Poisson transform was given as a fibre-integral in Definition $2.1.$ It appears from the perspective of the above functional equation that an explicit expression for each of the Fourier and Helgason Fourier transforms {\it is only possible} in the presence of an {\it invertible} Poisson transform.

\indent $4.$ There is the need to talk about the double fibre-integrals suggested by the composition $(\Phi_{k,l}\circ\widehat{f}_{l,k})(x,b)=\Phi_{k,l}(\widehat{f}_{l,k}(x,b))$ as regards the possibility of a {\it Fubini} argument. It is however known that the $G$ action on the product space $G/K\times G/P$ is transitive and that the kernel of the point $(eK,eP)$ is exactly $M:=K\cap P,$ leading to the diffeomorphism $G/M\widetilde{\longrightarrow} G/K\times G/P$ of $G-$manifolds, $[4.],\;p.\;46.$ This allows us to write (with appropriate {\it weighted} measures) the double fibre-integrals $\int_{G/P}\int_{G/K}$ as $\int_{G/M},$ and as $\int_{G/K}\int_{G/P}.$
 	
\ \\
	\indent {\bf \S4. Fourier and Helgason Fourier transform Maps for Invertible Poisson transforms.}
	
	\indent We now set out to invert the vector bundle-valued Poisson transform $\Phi_{k,l}$ on $G/P.$ However, since every differential form has a uniquely defined local expression (in terms of some local charts), we conclude that $$\Phi_{k,l}(\alpha)=\sum^{l}_{i=1}(P_{T}\alpha)_{_{i}}dy_{i_{1}}\wedge\cdots\wedge dy_{i_{l}},$$ where $\alpha\in\Omega^{k}(G/P,V[\chi]).$ This is a due to the fact that $dy_{i_{1}}\wedge\cdots\wedge dy_{i_{l}},$ $1\leq i_{1}<\cdots i_{l} \leq n,$ form a basis for $\Omega^{l}(G/K,W).$  Here $P_{T}$ is some $V-$valued smooth function on $\Omega^{k}(G/P,V[\chi])$ which we now study in details.
	
	\indent Consider an arbitrary irreducible unitary representation $(\tau,V)$ of $K.$ The space of the smooth sections of the associated vector bundle $G\times_{K}V$ shall be denoted by $\Gamma(G\times_{K}V)$ and may be identified with the space $C^{\infty}Ind^{G}_{K}(\tau)$ defined as $$C^{\infty}Ind^{G}_{K}(\tau)=\{f\in C^{\infty}(G/K,V):f(xk^{-1})=\tau(k)^{-1}f(x),\;\forall\;x\in G/K,k\in K\}.$$

\indent Let $$V=\bigoplus_{\sigma\in\widehat{M}}V(\sigma)$$ denote the decomposition into the $M-$isotopic parts, $V(\sigma),$ in which we write $\sigma\in \tau$ whenever $V(\sigma)\neq0.$ Let $\Lambda_{\sigma}$ denote the sum of the highest weight of $\sigma$ and $\rho_{o}.$
	
	\indent For any finite-dimensional representation $(\delta,V_{\delta})$ of the minimal parabolic subgroup $B=MAN$ of $G,$ we shall set $C^{\infty}Ind^{G}_{P}(\delta)$ to mean the vector space $C^{\infty}Ind^{G}_{P}(\delta)=$ $$\{\alpha\in C^{\infty}(G/K,V):\alpha(xman)=e^{-\rho\log a}\delta^{-1}(man)\alpha(x),\;\forall x\in G/K,\;man\in P\}$$ (endowed with the topology induced from $C^{\infty}(G/K,V_{\delta})$) and define a Poisson transform as a continuous linear $G-$equivariant map from $C^{\infty}Ind^{G}_{P}(\delta)$ into $C^{\infty}Ind^{G}_{K}(\tau).$ It can be shown that every $T\in Hom_{M}(V_{\delta},V)(=Hom(V_{M},W_{M}))$ corresponds to a Poisson transform $P_{T}$ which is well-defined on every $\alpha\in C^{\infty}Ind^{G}_{P}(\delta)$ as $$P_{T}\alpha(gK)=\int_{K}\tau(k)T(\alpha(gk))dk=\int_{K}(\tau(k)\circ T)(\alpha(gk))dk,$$ where $dk$ is the normalized Har measure on the compact Lie group $K,\;[7.]$ and that $T\mapsto P_{T}$ is a bijection from $Hom_{M}(V_{\delta},V)$ {\it onto} the $V-$valued space $\mathfrak{P}_{Hom_{M}(V_{\delta},V)}(C^{\infty}Ind^{G}_{P}(\delta))$ of Poisson transforms of vector-valued functions, $[12.].$ We can extend this result to include the vector bundle-valued case.
	
	\indent {\bf 4.1 Proposition.} Denote by $\mathfrak{P}_{Hom(V_{M},W_{M})}(\Omega^{k}(G/P,V[\chi]))$ the space of all vector bundle-valued Poisson transforms of $k-$forms on $G/P$ corresponding to the kernel $\phi_{k,l}\in\Omega^{l,n-k}(G/K/\times G/P,Hom(V_{M},W_{M}))$ as given in Definition $2.1.$ Then the map $\phi_{k,l}\mapsto\Phi_{k,l}$ is a bijection from the space $\Omega^{l,n-k}(G/K/\times G/P,Hom(V_{M},W_{M}))$ onto $\mathfrak{P}_{Hom(V_{M},W_{M})}(\Omega^{k}(G/P,V[\chi])).$ In particular, $\mathfrak{P}_{Hom_{M}(V_{\delta},V)}(C^{\infty}Ind^{G}_{P}(\delta))$ is a closed subspace of the vector space $\mathfrak{P}_{Hom(V_{M},W_{M})}(\Omega^{k}(G/P,V[\chi])).$
	
	\indent {\bf Proof.} Note that the vector space $\mathfrak{P}_{Hom(V_{M},W_{M})}(\Omega^{k}(G/P,V[\chi]))$ is explicitly given as $\mathfrak{P}_{Hom(V_{M},W_{M})}(\Omega^{k}(G/P,V[\chi]))=\{\Phi_{k,l}(\alpha)\in\Omega^{l}(G/K,W):\Phi_{k,l}(\alpha)=\sum^{l}_{i=1}(P_{T}\alpha)_{_{i}}dy_{i_{1}}\wedge\cdots\wedge dy_{i_{l}},\alpha\in C^{\infty}Ind^{G}_{P}(\delta),\;T\in Hom_{M}(V_{\delta},V)\}.$ Here each $T\in Hom_{M}(V_{\delta},V)$ is the $0-$forms among the members $\phi_{k,l}$ in $\Omega^{l,n-k}(G/K/\times G/P,Hom(V_{M},W_{M})).$ Hence the vector space denoted as $\mathfrak{P}_{Hom_{M}(V_{\delta},V)}(C^{\infty}Ind^{G}_{P}(\delta))$ consists of all the $0-$forms among all of the $l-$forms contained in the vector space $\mathfrak{P}_{Hom(V_{M},W_{M})}(\Omega^{k}(G/P,V[\chi])).$
	
	\indent Consider now the $0-$to-$(l,n-k)$ formation of differential-form map $$\mu:Hom_{M}(V_{\delta},V)\rightarrow Hom(V_{M},W_{M})$$ defined as $$\mu(\phi(\cdot))=\sum^{k}_{i=1}\phi(\cdot)_{_{i}}dx_{i_{1}}\wedge\cdots\wedge dx_{i_{k}}$$ which extends to the map $$\bar{\mu}:\mathfrak{P}_{Hom_{M}(V_{\delta},V)}(C^{\infty}Ind^{G}_{P}(\delta))\rightarrow \mathfrak{P}_{Hom(V_{M},W_{M})}(\Omega^{k}(G/P,V[\chi]))$$ given as $$\bar{\mu}(P_{T}(\cdot))=\sum^{k}_{i=1}P_{T}(\cdot)_{_{i}}dx_{i_{1}}\wedge\cdots\wedge dx_{i_{k}}.$$ Both $\mu$ and $\bar{\mu}$ are bijections, which is due to the uniqueness of the formation of local expressions for differential forms. We set $$\kappa:Hom_{M}(V_{\delta},V)\longrightarrow \mathfrak{P}_{Hom_{M}(V_{\delta},V)}(C^{\infty}Ind^{G}_{P}(\delta))$$ as well as $$\kappa^{*}:Hom(V_{M},W_{M})\longrightarrow \mathfrak{P}_{Hom(V_{M},W_{M})}(\Omega^{k}(G/P,V[\chi])).$$ It is known, from $[12.],$ p. $863,$ that $\kappa$ is a bijection. Since $\kappa^{*}=\mu\circ\kappa\circ\mu^{-1}$ we have the result$.\;\Box$
	
	\indent The bijectivity of the map $\bar{\mu}$ implies that an understanding of the vector space $\mathfrak{P}_{Hom(V_{M},W_{M})}(\Omega^{k}(G/P,V[\chi]))$ (of $V-$valued $l-$forms) completely follows from an understanding of the vector space $\mathfrak{P}_{Hom_{M}(V_{\delta},V)}(C^{\infty}Ind^{G}_{P}(\delta)),$ which we now focus on.
	
	\indent Let $g\in G.$ Denote by $||g||$ the operator norm of Ad$(g)$ on $\mathfrak{g}$ and denote $\langle X,Y\rangle_{\theta}$ by the {\it Cartan inner product} $$\langle X,Y\rangle_{\theta}=-k(X,\theta Y),$$ for the Cartan-Killing form $k$ on $\mathfrak{g}.$ We shall refer to a $V-$valued fnction $f$ as {\it increasing at most exponentially} if $||f||_{r}<\infty,$ where $$||f||_{r}=\sup_{g\in G}||g||^{-r}|f(g)|,$$ $r\in\mathbb{R}.$ Now set $$C_{r}(G,V)=\{f\in C(G,V): ||f||_{r}<\infty\},$$ for some, $r\in\mathbb{R}.$ It is known that for each $\alpha\in C^{\infty}Ind^{G}_{P}(\sigma\otimes(-\lambda)\otimes1)$ and $T\in Hom_{M}(V_{\delta},V),$ we always have that $P_{T}\alpha\in C_{r(\lambda)}(G,V),$ where $$r(\lambda)=C|\Re\lambda-\rho|,$$ for some constant $C.$
	
	\indent Let $\mathcal{E}_{\lambda-\Lambda}Ind^{G}_{K}(\tau)$ denote the eigenspace of $Z_{\tau}$ (whose character is $\lambda-\Lambda,$ where $\lambda\in\mathfrak{a}^{*}_{\mathbb{C}}$ and $\Lambda$ is the infinitesimal character of an irreducible representation of $M$ in $\tau$). Set $$\mathcal{E}^{\infty}_{\lambda-\Lambda}Ind^{G}_{K}(\tau)=\bigcup_{r\in\mathbb{R}}\mathcal{E}^{\infty}_{\lambda-\Lambda,r}Ind^{G}_{K}(\tau),$$ where $\mathcal{E}^{\infty}_{\lambda-\Lambda,r}Ind^{G}_{K}(\tau)=\mathcal{E}_{\lambda-\Lambda}Ind^{G}_{K}(\tau)\bigcap C^{\infty}_{r}(G,V).$ We now state a major asymptotic result for members of $\mathcal{E}_{\lambda-\Lambda}Ind^{G}_{K}(\tau).$
	
	\indent {\bf 4.2 Theorem.} ($[12.],$ p. $870$). Define the set $$X(\lambda,\Lambda)=\{\omega(\lambda-\Lambda)+\Lambda_{\sigma}-\rho-\mathbb{N}\cdot\triangle: \omega\in\widetilde{W},\sigma\in\tau\;\mbox{and}\;(\omega(\lambda-\Lambda)+\Lambda_{\sigma})_{_{|_{\mathfrak{t}}}}=0\}.$$
	
	\indent $(i)$ For each $f\in \mathcal{E}^{\infty}_{\lambda-\Lambda}Ind^{G}_{K}(\tau),$ $g\in G$ and $\xi\in X(\lambda,\Lambda),$ there exists a unique $V-$valued polynomial $p_{\lambda,\xi}(f,g,\cdot)$ on $\mathfrak{a}$ such that $$f(tH)\sim\sum_{\xi\in X(\lambda,\Lambda)}p_{\lambda,\xi}(f,g,tH)e^{t\xi(H)},\;\;t\rightarrow\infty$$ at every $H\in\mathfrak{a}^{+},$ in which the polynomials are of degree  $\leq d,$ $d$ being the number of elements in $\sum^{+}(\mathfrak{g}_{\mathbb{C}},\mathfrak{h}_{\mathbb{C}}).$
	
	\indent $(ii)$ If $r\in\mathbb{R}$ and $\xi\in X(\lambda,\Lambda)$ there exists $r^{\prime}\in\mathbb{R}$ such that $f\mapsto p_{\lambda,\xi}(f,\cdot,\cdot)$ is a continuous map of $\mathcal{E}^{\infty}_{\lambda-\Lambda}Ind^{G}_{K}(\tau)$ into $C^{\infty}_{r}(G,V)\otimes P_{d}(\mathfrak{a}),$ equivariant for the left action of $G$ on $\mathcal{E}^{\infty}_{\lambda-\Lambda}Ind^{G}_{K}(\tau)$ to $C^{\infty}_{r}(G,V)\otimes P_{d}(\mathfrak{a}).\;\;\Box$
	
	\indent It is known that $p_{\lambda,\xi}(f,g,\cdot)$ is a constant on $\mathfrak{a}$ for each $g\in G,$ which we then simply denote by $p_{\lambda,\xi}(f,\cdot).$ Now let $\alpha^{V}=\frac{2\alpha}{\langle\alpha,\alpha\rangle}$ and consider the spaces $$\mathcal{U}_{1}=\{\lambda-\Lambda:\lambda\in\mathfrak{a}^{*}_{\mathbb{C}},\Lambda\in\mathfrak{t}^{*}_{\mathbb{C}}, \langle\lambda-\Lambda,\alpha^{V}\rangle\notin\mathbb{Z},\forall\;\alpha\sum(\mathfrak{g}_{\mathbb{C}},\mathfrak{h}_{\mathbb{C}}),\alpha_{_{|_{\mathfrak{a}}}}\ne0\}$$ and $$\mathcal{U}_{2}=\{\lambda\in\mathfrak{a}^{*}_{\mathbb{C}}:\langle \alpha,\gamma^{V} \rangle\notin-\mathbb{N},\;\forall\;\gamma\in\sum(\mathfrak{g},\mathfrak{a})\},$$ and the {\it boundary value map} $\beta_{\lambda}$ on $\mathcal{E}^{\infty}_{\lambda-\Lambda}Ind^{G}_{K}(\tau),$ given by $$\beta_{\lambda}(f)=p_{\lambda,\lambda-\rho}(f,\cdot),$$ $f\in \mathcal{E}^{\infty}_{\lambda-\Lambda}Ind^{G}_{K}(\tau).$
	
	\indent If $\tau(\Lambda)$ denote the restriction of $\tau$ to $M$ with representation space $V(\Lambda)$ then $\beta_{\lambda}$ is a linear, continuous and $G-$equivariant map of $\mathcal{E}^{\infty}_{\lambda-\Lambda}Ind^{G}_{K}(\tau)$ into $C^{\infty}Ind^{G}_{P}(\tau(\Lambda)\otimes(-\lambda)\otimes1),$ for each $r\in \mathbb{R}$ and if $\Omega\in\mathfrak{a}^{*}_{\mathbb{C}}$ is open with an holomorphic family $\{f_{\lambda}\}_{\lambda\in\Omega}$ in $\mathcal{E}^{\infty}_{\lambda-\Lambda}Ind^{G}_{K}(\tau),$ then $\lambda\rightarrow\beta_{\lambda}(f_{\lambda})$ is holomorphic in $\Omega\bigcap\mathcal{U}_{2}.$ We also have that $$\lim_{t\rightarrow\infty}e^{(-\lambda+\rho)(tH)}f(g\exp tH)=(\beta_{\lambda }f)(g),$$ whenever $\lambda-\Lambda\in \mathcal{U}_{1},$ $\Re\langle\lambda,\alpha\rangle\textgreater0$ for each $\alpha\in\sum^{+}(\mathfrak{g},\mathfrak{a}),$ $H\in\mathfrak{a}^{+}$ and $f\in\mathcal{E}^{\infty}_{\lambda-\Lambda}Ind^{G}_{K}(\tau).$ The major result on the inversion of the Poisson map $\alpha\mapsto P_{\lambda}\alpha$ where $\alpha\in C^{\infty}Ind^{G}_{P}(\tau(\Lambda)\otimes(-\lambda)\otimes1)$ and $$(P_{\lambda}\alpha)(gK)=\int_{K}\tau(k)T_{\lambda}\alpha(gk)dk$$ is the following.
	
	\indent {\bf 4.3 Theorem.} ($[12.],$ p. $879$). Let $\lambda-\Lambda\in \mathcal{U}_{1}$ and $\lambda\in \mathcal{U}_{2}.$ Let $C_{o}(\lambda)$ denote the restriction to $V(\Lambda)$ of the {\it generalized Harish-Chandra $C-$function,} where the map $\lambda\rightarrow C(\lambda)$ is given by $C(\lambda)=\int_{\Theta(N)}e^{-(\lambda+\rho)(H(\bar{n}))}\tau(k(\bar{n}))d\bar{n},$ and let $$P_{\lambda}:C^{\infty}Ind^{G}_{P}(\tau(\Lambda)\otimes(-\lambda)\otimes1)\rightarrow\mathcal{E}^{\infty}_{\lambda-\Lambda}Ind^{G}_{K}(\tau)$$ be as defined above. Then,
	
	\indent $(i)$ $\beta_{\lambda}P_{\lambda}\alpha=C_{o}(\lambda)\alpha,$ for each $\alpha\in C^{\infty}Ind^{G}_{P}(\tau(\Lambda)\otimes(-\lambda)\otimes1);$
	
	\indent $(ii)$ $\beta_{\lambda}$ (respectively, $P_{\lambda}$) is surjective (respectively, injective) whenever $\det C_{o}(\lambda)\neq0.$ Indeed, $P_{\lambda}$ is a bijection and its inverse, $P^{-1}_{\lambda},$ is given as $C_{o}(\lambda)^{-1}\circ\beta_{\lambda}.\;\;\Box$
	
	\indent In this context, the Helgason Fourier transform of every vector bundle-valued function $f\in C^{\infty}Ind^{G}_{P}(\tau(\Lambda)\otimes(-\lambda)\otimes1)$ is now seen as $$\widehat{f}_{l,k,\varepsilon(\lambda)}(x,b)=P^{-1}_{\lambda}\circ(\varphi^{U^{\sigma\nu}}_{\lambda,l,k}*f)(x)=C_{o}(\lambda)^{-1}\circ\beta_{\lambda}\circ(\varphi^{U^{\sigma\nu}}_{\lambda,l,k}*f)(x)$$ while the Fourier transform is then $$f^{\triangle}_{l,k,\upsilon(\lambda)}(b,x)=\int_{K}\tau(k)\circ T_{\lambda}\circ\widehat{f}_{l,k,\varepsilon(\lambda)}(gk,b)dk$$ $$\;\;\;\;\;\;\;\;\;\;\;\;\;\;\;\;\;\;\;\;\;\;\;\;\;\;\;\;\;\;\;\;\;\;\;\;\;\;\;\;\;\;\;=\int_{K}\tau(k)\circ T_{\lambda}\circ C_{o}(\lambda)^{-1}\circ\beta_{\lambda}\circ(\varphi^{U^{\sigma\nu}}_{\lambda,l,k}*f)(gk)dk$$
	
	\indent We now return to the vector space $\mathfrak{P}_{Hom(V_{M},W_{M})}(\Omega^{k}(G/P,V[\chi]))$ consisting of all $W-$valued $l-$forms Poisson transform image of $V-$valued $k-$forms in $\Omega^{k}(G/P,V[\chi])$ and the bijective map $\bar{\mu}$ of Proposition $4.1.$ Define the vector bundle-valued Harish-Chandra $c-$functions $${\bf C_{o}}(\lambda)={\bf C}(\lambda)_{_{|_{V(\Lambda)}}},$$ where $${\bf C}(\lambda)=\sum^{k}_{i=1}C(\lambda)_{_{i}}dx_{i_{1}}\wedge\cdots\wedge dx_{i_{k}}$$ and $$\beta^{V}_{\lambda}(\omega)=\sum^{k}_{i=1}\beta_{\lambda}(\omega)_{_{i}}dx_{i_{1}}\wedge\cdots\wedge dx_{i_{k}},$$ where $$\omega=\sum^{k}_{i=1}f_{_{i}}dx_{i_{1}}\wedge\cdots\wedge dx_{i_{k}}$$ for any local chart. It is known that this representation by a local chart is unique, independent of the choice of chart and that the wedge products, $dx_{i_{1}}\wedge\cdots\wedge dx_{i_{k}},$ $1\leq i_{1}<\cdots< i_{k}\leq k,$ form a basis for $\Omega^{k}(G/P,V[\chi]).$

\indent We also set $$\mathcal{E}^{\infty,V}_{\lambda-\Lambda}Ind^{G}_{K}(\tau)=\bigcup_{r\in\mathbb{R}}\mathcal{E}^{\infty,V}_{\lambda-\Lambda,r}Ind^{G}_{K}(\tau),$$ in which we have $\mathcal{E}^{\infty,V}_{\lambda-\Lambda,r}Ind^{G}_{K}(\tau)=\{\alpha\in\Omega^{k}(G/P, V[\chi]): \alpha=\sum^{k}_{i=1}f_{_{i}}dx_{i_{1}}\wedge\cdots\wedge dx_{i_{k}}\;\mbox{and}\;f\in\bigcup_{r\in\mathbb{R}}\mathcal{E}^{\infty}_{\lambda-\Lambda,r}Ind^{G}_{K}(\tau)\}.$ If $\det{\bf C_{o}}(\lambda)\neq0,$ we shall also set $${\bf C_{o}}(\lambda)^{-1}=\sum^{k}_{i=1}C(\lambda)^{-1}_{_{i}}dy_{i_{1}}\wedge\cdots\wedge dy_{i_{k}}$$ such that $$\sum^{k}_{i=1}dx_{i_{1}}\wedge\cdots\wedge dx_{i_{k}}\wedge(dy_{i_{1}}\wedge\cdots\wedge dy_{i_{k}})=\sum^{k}_{i=1}\delta_{i_{k}j_{k}}=1.$$
	
	\indent Now let $\Phi_{k,l,\lambda}$ denote the Poisson transform corresponding to some $\lambda\in\mathfrak{a}^{*}_{\mathbb{C}}.$ An analogous result to Theorem $4.3$ is the following.
	
	\indent {\bf 4.4 Theorem.} Let $\lambda-\Lambda\in \mathcal{U}_{1}$ and $\lambda\in \mathcal{U}_{2}.$ Then,
	
	\indent $(i)$ $\beta^{V}_{\lambda} \Phi_{k,l,\lambda}(\alpha)={\bf C_{o}}(\lambda)\alpha,$ for each $\alpha\in \Omega^{k}(G/P,V[\chi]);$
	
	\indent $(ii)$ $\beta^{V}_{\lambda}$ (respectively, $\Phi_{k,l,\lambda}(\alpha)$) is surjective (respectively, injective) whenever $\det {\bf C_{o}}(\lambda)\neq0.$ Indeed, $\Phi_{k,l,\lambda}(\alpha)$ is a bijection and its inverse, $\Phi^{-1}_{k,l,\lambda}(\alpha),$ is given as ${\bf C_{o}}(\lambda)^{-1}\circ\beta^{V}_{\lambda}.$
	
	\indent {\bf Proof.} The bijectivity in the closed diagram in the proof of Proposition $4.1$ proves that Theorem $4.3 (i)$ (respectively, Theorem $4.3(ii)$) implies $(i)$ (respectively, $(ii)$) of the present result$.\;\;\Box$
	
	\indent The import of these results, especially Theorem $4.4$ on the Helgason Fourier transform is being able to get an explicit expression for this transform from $$\widehat{f}_{l,k}(x,b)=\Phi^{-1}_{k,l}((\varphi^{U^{\sigma\nu}}_{\tau,l,k}*f)(x)).$$
	
	\indent {\bf 4.5 Proposition.} Let $\lambda-\Lambda\in \mathcal{U}_{1}$ and $\lambda\in \mathcal{U}_{2}.$ The Helgason Fourier transform $f\mapsto\widehat{f}_{l,k,\varepsilon(\lambda)}$ is explicitly given as $$\widehat{f}_{l,k,\varepsilon(\lambda)}(x,b)=({\bf C_{o}(\lambda)}^{-1}\circ\beta^{V}(\lambda))\circ((\varphi^{U^{\sigma\nu}}_{\lambda,l,k}*f)(x)),$$ $(x,b)\in G/K\times G/P,$ for some $\lambda-$linear relation $\varepsilon(\lambda).\;\;\Box$
	
	\indent In other words, if $\lambda-\Lambda\in \mathcal{U}_{1}$ and $\lambda\in \mathcal{U}_{2},$ then $$\widehat{f}_{l,k,\varepsilon(\lambda)}(x,b)=({\bf C_{o}(\lambda)}^{-1}\circ\beta^{V}(\lambda))\circ(\int_{G/K}\varphi^{U^{\sigma\nu},t}_{\lambda,l,k}\wedge\pi^{*}_{K}f)(x),$$ $(x,b)\in G/K\times G/P,$ for some $\lambda-$linear relation $\varepsilon(\lambda).$
\ \\
\ \\
\ \\
	\indent {\bf 4.6 Remarks.}
	
	\indent $(1)$ The above proposition gives the Helgason Fourier transform of $W-$ valued $l-$forms on $G/K$ as a fibre integral. Further explication of this transform beyond Proposition $4.5$ would require an explicit expression for the convolution $$\varphi^{U^{\sigma\nu}}_{\lambda,l,k}*f:=\int_{G/K}\varphi^{U^{\sigma\nu},t}_{\lambda,l,k}\wedge\pi^{*}_{K}f,$$ which in turn requires an explicit formula for the translated $W-$valued $\tau-$spherical differential form $\varphi^{U^{\sigma\nu},t}_{\lambda,l,k}$ on $G/K,$ as well as the fibre integral formulae for both ${\bf C_{o}(\lambda)}^{-1}$ and $\beta^{V}(\lambda)$ and the $\lambda-$linear relation $\varepsilon(\lambda).$
	
	\indent $(2)$ The classical examples of the Fourier, Helgason Fourier and Poisson transforms are trivial cases in the labyrinth of the theory constructed for the three transforms so far. Indeed, it was pointed out in $[4.]$ that one of the cases that could reduce to the classical Poisson transform is from a uniformly constructed Poisson kernel $\phi_{k}$ of degree $(k,n-k)$ given as $$\phi_{k}=*_{K}(\tau\wedge(d_{P}E^{*})^{n-k}),$$ with $k=0,\cdots,n=\dim(G/P),$ where $\tau\in\Lambda^{d}((\mathfrak{g}/\mathfrak{m})^{*}_{o})$ is an $M-$invariant element extended trivially to an invariant $d-$form on $\mathfrak{g}/\mathfrak{m}$ and $E^{*}\in(\mathfrak{g}/\mathfrak{m})^{*}$ is an induced non-trivial $\mathfrak{m}-$invariant element of degree $(1,0)$ with trivial $K-$derivative but with $P-$derivative known as $$d_{P}E^{*}(F_{X},G_{Y})=i(\sum^{k}_{j=-k}j^{2}\dim(\mathfrak{g}_{j}))^{-1}k(X,Y),$$ $X\in \mathfrak{g}_{i}, Y\in \mathfrak{g}_{-i}.$ If $\Phi_{k}$ denote the Poisson transform corresponding to the kernel $\phi_{k},$ it is noted that $\Phi_{k}$ is {\it coclosed,} $[4.],$ p. $61.$
	
	\indent With $k=0,$ the classical Poisson transform, $P_{\lambda^{\prime}},$ $\lambda^{\prime}\in\mathfrak{a}^{*}_{\mathbb{C}},$ is given as $$P_{-i\lambda^{\prime}}=\Phi_{o}\circ \chi_{\rho-\lambda^{\prime}},$$ where $\Phi_{k=0}$ is the Poisson transform corresponding to the above kernel $\phi_{k=0}$ and $\chi_{\lambda^{\prime}}:G\rightarrow \mathbb{R}[\chi]$ is the smooth map given as $$\chi_{\lambda^{\prime}}(g)=e^{-{\lambda^{\prime}}(H(g))}\widetilde{f}(k(g)),$$ $\widetilde{f}(k)=f(kM).$ Indeed, $\triangle_{K}\Phi_{o}(\sigma)=-\langle\lambda^{\prime},\lambda^{\prime}-2\rho\rangle\Phi_{o}(\sigma)$ for any $\lambda^{\prime}-$density $\sigma.$ In particular, $$\triangle_{K}P_{-i\lambda^{\prime}}=\triangle_{K}\Phi_{0}(\chi_{\rho-\lambda^{\prime}})$$ $$\;\;\;\;\;\;\;\;\;\;\;\;\;\;\;\;\;\;\;\;\;\;\;\;\;\;\;\;\;\;=-\langle\lambda^{\prime},\lambda^{\prime}-2\rho\rangle\Phi_{o}(\chi_{\rho-\lambda^{\prime}})$$ $$\;\;\;\;\;\;\;\;\;\;\;\;\;\;\;\;\;\;\;\;\;\;\;\;=-\langle\lambda^{\prime},\lambda^{\prime}-2\rho\rangle P_{-i\lambda^{\prime}}.$$ By considering complex line bundles and the weight $\rho-i\lambda$ (for $\lambda\in\mathfrak{a}^{*}_{\mathbb{C}}$) as usually done in the classical theory of harmonic analysis, it follows that the above eigenvalue $-\langle\lambda^{\prime},\lambda^{\prime}-2\rho\rangle$ of $P_{-i\lambda^{\prime}}$ becomes (with $\lambda^{\prime}=\rho-i\lambda$), $$-\langle\lambda^{\prime},\lambda^{\prime}-2\rho\rangle=-\langle\rho-i\lambda,\rho-i\lambda-2\rho\rangle=-\langle\rho-i\lambda,-\rho-i\lambda\rangle$$ $$\;\;\;\;\;\;\;\;\;\;\;\;\;\;\;\;\;\;\;\;\;\;\;\;=-\langle\rho,-\rho\rangle-\langle\rho-i\lambda\rangle-\langle-i\lambda,-\rho\rangle-\langle-i\lambda,-i\lambda\rangle$$ $$\;\;\;\;\;=\langle\rho,\rho\rangle+i\langle\lambda,\rho\rangle-i\langle\lambda,\rho\rangle-i^{2}\langle\lambda,\lambda\rangle$$ $$=\langle\rho,\rho\rangle+\langle\lambda,\lambda\rangle,\;\;\;\;\;\;\;\;\;\;\;\;\;\;\;\;\;\;\;\;\;\;\;\;\;$$ as contained in \S 2. above and in $[5.].$
	
	\indent $(3)$ The present paper {\it weans} the theory of Fourier and Helgason Fourier transforms from the classical perspective of being considered as some individually defined  transforms constructed in the fashion of the past cases. Presently we have the {\it abstract,} the {\it coclosed} and the {\it Berstein-Gelfand-Gelfand (BGG)-compatible} types to choose from. The appropriate outlook is to seek how much of the properties of the classical cases (and in the presence of what conditions) would be valid for any of these types.
	
	\indent {\bf 4.7 Further Directions.}
	
	\indent $(1)$ A pertinent question may then be asked; What are the distinguishing features of the Helgason Fourier transform $\widehat{f}$ among all of the members of $\Omega^{k}(G/K\times G/P,V[\chi])?$ A first answer is that: only the $\alpha=\widehat{f}$ satisfying Definition $3.2.$ That is, the functional equation $$(\varphi^{U^{\sigma\nu}}_{\tau,l,k}*f)(x)=\Phi_{k,l}(\alpha)$$ solves only to give $\alpha=\widehat{f}.$ Thus there may be other members of $\Omega^{k}(G/K\times G/P,V[\chi])$ for which there is no single $f\in\Omega^{l}(G/K,W)$ satisfying the functional equation. Hence we could consider the subspace $\Omega^{k}_{H}(G/K\times G/P,V[\chi])$ of $\Omega^{k}(G/K\times G/P,V[\chi])$ consisting of $V[\chi]-$valued differential forms on $G/K\times G/P$ that are Helgason Fourier transforms of some members of $\Omega^{l}(G/K,W).$ This brings up the vector subspace $\Omega^{k}_{H}(G/K\times G/P,V[\chi])=\{\alpha\in\Omega^{k}_{H}(G/K\times G/P,V[\chi]):\alpha=\widehat{f},\;\mbox{for some}\;f\in\Omega^{l}(G/K,W)\}.$ Clearly, $$\Omega^{k}_{H}(G/K\times G/P,V[\chi])\subseteq\Omega^{k}(G/K\times G/P,V[\chi])\subseteq\Omega^{k}(G/P,V[\chi])$$ so that $\Phi_{k,l}$ could be restricted to $\Omega^{k}_{H}(G/K\times G/P,V[\chi]),$ which now carries the induced topology from $\Omega^{k}(G/K\times G/P,V[\chi]).$
	
	\indent $(2)$ The present consideration may be used to consider the Fourier, Poisson and Helgason Fourier transforms in the context of vector bundled-valued harmonic analysis of {\it spherical functions of arbitrary type associated to a representation $\pi$} as follows. Let $\delta\in\widehat{K}$ whose character is denoted as $\chi_{\delta}$ and define $$\xi_{\delta}(k)=\dim(\delta)\chi_{\delta}(k^{-1}),\;k\in K.$$ Choose any representation $\pi$ of $K$ in a (complete) Hausdorff locally convex space $V$ and set $E_{\delta}=\pi(\chi_{\delta}).$ Each $E_{\delta}$ is a continuous projection and for any finte set $F\subset\widehat{K},$ we may define $E_{F}=\sum_{\delta\in F}E_{\delta}$ and set $V_{F}=E_{F}(V).$
	
	\indent The spherical function of type $F$ associated to $\pi,$ written as $\varphi_{\pi,F}$ is given by $$\varphi_{\pi,F}(x)=E_{F}\pi(x)E_{F},\;x\in G.$$ We shall set $\widetilde{f}=f(x^{-1})^{conj},$ $x\in G$ and denote the right regular representation of $G$ by $(r(x)f)(y)=f(yx),$ $x,y\in G.$ The {\it Fourier transform on $G$ of arbitrary type $F$ associated to $\pi,$} written as $f\mapsto f^{\triangle}(\pi,x),$ may in the present case be given as $$f^{\triangle}(\pi,x)=(\varphi_{\pi,F}*f)(x),$$ $x\in G.$ The rich structure of this transform is immediately seen in the fact that (for $f\in C_{c}(G)$) $$f^{\triangle}(\pi,x)=(\varphi_{\pi,F}*f)(x)=\int_{G}\varphi_{\pi,F}(y)f(y^{-1}x)dy$$ $$=\int_{G}\varphi_{\pi,F}(y)(r(x)f)(y^{-1})dy=\int_{G}\varphi_{\pi,F}(y)\widetilde{(r(x)f)(y)^{conj}}dy$$ $$=\int_{G}E_{F}\pi(x)E_{F}\widetilde{(r(x)f)(y)^{conj}}dy=\pi(r(x)f^{conj})_{_{|_{_{V_{F}}}}}=\pi(r(x)f^{conj})^{\dagger}_{_{|_{_{V_{F}}}}},$$ where $r(x)f^{conj}\in\xi_{F}*C_{c}(G)*\xi_{F}$ and $\dagger$ represents the Hilbert space adjoint. This observation opens up the subject of Fourier transform of types $F$ to the consideration of Poisson and Helgason Fourier transforms of type $F$ associated to the representation $\pi.$
	
    \indent $(3)$ Beyond the Helgason Fourier transform theory, the functional equation of this paper may be considered in a more general form as seeking all those $\alpha\in\Omega^{k}(G/K\times G/P,V)$ for which $$(\varphi^{U^{\sigma\nu}}_{\tau,l,k}*f)(x)=\Phi_{k,l}(\alpha(x,b)),$$ where $(x,b)\in G/K\times G/P$ for some $f\in\Omega^{l}(G/K,W)$ (or the inverse problem of seeking all those $f\in\Omega^{l}(G/K,W)$ satisfying the same equation for some $\alpha\in\Omega^{k}(G/K\times G/P,V)$). We could then have the general transform-map $$f\mapsto\alpha:\Omega^{l}(G/K,W)\rightarrow\Omega^{k}(G/K\times G/P,V)$$ which should subsume the theory of Helgason Fourier transform.
    
    \indent In complete analogy with Proposition $4.5,$ we equally have an explicit fibre-integral expression for the Fourier transform.
    
    \indent {\bf 4.8 Theorem.} Let $\lambda-\Lambda\in \mathcal{U}_{1}$ and $\lambda\in \mathcal{U}_{2}.$ Then the Fourier transform $f\mapsto f^{\triangle}_{l,k,\upsilon(\lambda)}$ is explicitly given as the double fibre-integrals $$f^{\triangle}_{l,k,\upsilon(\lambda)}(b,x)=\int_{G/P}\phi_{k,l,\lambda}\wedge\pi^{*}_{P}(({\bf C_{o}(\lambda)}^{-1}\circ\beta^{V}(\lambda))\circ(\int_{G/K}\varphi^{U^{\sigma\nu},t}_{\lambda,l,k}\wedge\pi^{*}_{K}f)(x)),$$ for all $(x,b)\in G/K\times G/P,$ for some $\lambda-$linear relation $\upsilon(\lambda).$
    
    \indent {\bf Proof.} $f^{\triangle}_{l,k,\upsilon(\lambda)}(b,x)=\Phi_{k,l,\lambda}(\widehat{f}_{l,k,\varepsilon(\lambda)}(x,b))=\int_{G/P}\phi_{k,l,\lambda}\wedge\pi^{*}_{P}\widehat{f}_{l,k,\varepsilon(\lambda)}(x,b)$ $$=\int_{G/P}\phi_{k,l,\lambda}\wedge\pi^{*}_{P}(({\bf C_{o}(\lambda)}^{-1}\circ\beta^{V}(\lambda))\circ(\int_{G/K}\varphi^{U^{\sigma\nu},t}_{\lambda,l,k}\wedge\pi^{*}_{K}f)(x)),$$ where $\varepsilon(\lambda)$ and $\upsilon(\lambda)$ are $\lambda-$linear relations$.\;\;\Box$
    
    \indent {\bf 4.9 Remarks.}
    
    \indent $(1)$ The above double fibre-integral formula for the Fourier transform is a confirmation and far reaching generalization of the double integrals of the same transform for the classical examples of $\mathbb{R},$ $\mathbb{R}^{n}$ and scalar-valued functions on $G/K$ earlier considered in the Introduction. This double fibre-integral formula was however not derived in the usual way of the past but via the functional equation of Definition $3.3.$ This approach has the benefits of reproducing all of the component terms of the transform, without having to impose any extra condition on the theory. Indeed, the fibre-integral expressions here derived for the Fourier and Helgason Fourier transforms give all of the component terms of these transforms and their strategic positions.
    
    \indent For the example of the Fourier transform in Theorem $4.8$ above, the component terms include the fibre-integral, $\int_{G/P},$ over $G/P;$ the corresponding Poisson kernel, $\phi_{k,l,\lambda};$ a wedge product; the pullback, $\pi^{*}_{P},$ of the projection $\pi_{P}:G/M\rightarrow G/P;$ inverse, ${\bf C_{o}(\lambda)}^{-1},$ of the invertible Harish-Chandra $C-$function, $\lambda\mapsto {\bf C_{o}}(\lambda);$ a composition; the boundary map $\beta^{V}_{\lambda};$ a composition; a fibre-integral, $\int_{G/K},$ over $G/K$ (evaluated at $x\in G/K$); the translation, $\varphi^{U^{\sigma\nu},t}_{\lambda,l,k},$ of the $Hom(V_{M},W_{M})-$ valued $\tau-$spherical $(m-l,k)-$form, $\varphi^{U^{\sigma\nu}}_{\lambda,l,k},$ on $G/K\times G/P;$ a wedge product and the image $\pi^{*}_{K}f$ of $f\in \Omega^{l}(G/K,W)$ under the pullback $\pi^{*}_{K}$ of the projection $\pi_{K}:G/M\rightarrow G/K.$ Apart from being able to now study the general properties of these transforms, these explicit expressions would clearly show which component term leads to which property of the Fourier or Helgason Fourier transform, we will be able to derive the condition(s) under which the known classical properties of these transforms are valid and to extend the theory to an arbitrary  Gelfand pair $(G,K).$
    
    \indent $(2)$ Theorem $3.7$ of $[8.]$ suggests that a Paley-Wiener theorem for the Fourier transform is essentially an injective Poisson-transform-image of the Paley-Wiener theorem for the Helgason Fourier transform. A programme for the vector bundle-valued $L^{2}-$harmonic analysis of the Fourier transform is also outlined in $[9.]$ and $[10.].$ Remarks $4.6(2)$ gives a uniform construction of the Poisson kernel, $\phi_{k},$ its corresponding Poisson transform, $\Phi_{k},$ and the scalar-valued Poisson transform, $$P_{-i\lambda^{\prime}}=\Phi_{o}\circ \chi_{\rho-\lambda^{\prime}},$$ that leads Proposition $4.5$ and Theorem $4.8$ back to the classical cases earlier outlined in the Introduction. It is now clear that {\it the Fourier transform is the Poisson-transform-completion of the Helgason Fourier transform.}

	\ \\
	\indent {\bf \S5. Conclusion.} The contribution of homological algebra (and its ramifications) to the {\it real analysis} of the Fourier and Helgason Fourier transforms is now accessible through the theory expounded in this paper, with attendant consequences on the classical cases of $\mathbb{R},$ $\mathbb{R}^{n}$ and the numerous Riemannian symmetric spaces $G/K.$ This is another instance where real, functional and harmonic analysis merge to be subsumed under the topics of homological algebra, as being sought in $[6.].$ In this regard, further homological properties of the component terms of both transforms would surely enrich the theory far beyond what classical analysis is capable of, into an aspect of mathematics that may be termed {\it Functional Homological Analysis.} The numerous applications of these transforms in diverse fields suggest the importance of such an undertaking.

	\ \\
	\indent {\bf \S6.  Declarations.}

	\indent {\it Ethics Approval and consent to participate:} I declare that there is no unethical conduct in the research and I also give my consent to participate in Ethics Approval.
	\ \\
	
	{\it Consent for publication:} I give my consent to the publication of this manuscript.
	\ \\
	
	{\it Availability of data and materials:} Materials consulted during the research are as listed in the References below.
	\ \\
	
	{\it Competing interest:} I declare that there is no competing interest in the conduct of this research.
	\ \\
	
	{\it Funding:} I declare that I received no funding for this research.
	\ \\
	
	{\it Authors' contributions:} I declare that I was the only one who contributed to the conduct of this research.
	\ \\
	
	{\it Acknowledgments:} I acknowledge authors whose research works were consulted as listed in the References below.

	\ \\
	\indent {\bf \S7. References.}
	\begin{description}
		
		\item [{[1.]}] Camporesi, R., The Helgason Fourier transform for homogeneous vector bundles over Riemannian symmetric spaces, {\it Pacific Journal of Mathematics,} vol. {\bf 179}, no. {\bf 2} $1997,$ p. $263-300.$
		
		\item [{[2.]}] Camporesi, R., The Helgason Fourier transform for homogeneous vector bundles over compact Riemannian symmetric spaces-the local theory, {\it J. Funct. Anal.,} vol. {\bf 220}, $2005,$ p. $97-117.$
		
		\item [{[3.]}] Chaib, A. O., On the Poisson transform on a homogeneous vector bundle over the quaternionic hyperbolic space, {\it Comptes Rendus Math$\acute{e}$matique,} vol. {\bf362}, $2024,$ p. $265-273.$
		
		\item [{[4.]}] Harrach, C., Poisson transform for differential forms adapted to the flat parabolic geometries on spheres, Ph. D thesis, University of Vienna, $2017.$
		
		\item [{[5.]}] Helgason, S., \textit{Geometric analysis on symmetric spaces,} Mathematical Survey and Monograph, AMS, $1994.$
		
		\item [{[6.]}] Hilgert, J., Poisson transforms and homological algebra, {\it Workshop Geometrie und Darstellungstheorie,} Beilefeld, $2022.$
		
		\item [{[7.]}] Olbrich, M., {\it Die Poisson-transformation f$\ddot{u}$r homogene Vektorb$\ddot{u}$ndel,} Ph. D thesis, Humboldt-Universit$\ddot{a}$t, Berlin, $1995.$
		
		\item [{[8.]}] Oyadare, O. O., A Paley-Wiener theorem for the joint-eigenspace Fourier transform on noncompact symmetric spaces.\\ Available at arXiv:$2409.09036$ [math.FA], $2024.$
		
		\item [{[9.]}] Oyadare, O. O., A note on the $L^{2}-$harmonic analysis of the joint-eigenspace Fourier transform.\\ Available at arXiv:$2410.09075$ [math.FA], $2024.$
		
		\item [{[10.]}] Oyadare, O. O., What is the JEFT?\\ Available at arXiv:$2412.10377$ [math.FA], $2024.$
		
		\item [{[11.]}] Sherman, T., The Helgason Fourier transform for compact Riemannian symmetric spaces of rank one, {\it Acta. Math.,} vol. {\bf 164}, $1990,$ p. $73-144.$
		
		\item [{[12.]}] Yang, A., Poisson transforms on vector bundles, {\it Trans. Amer. Math. Soc.,} vol. {\bf 350}, no. {\it 2}, $1998,$ p. $857-887.$

	\end{description}
\end{document}